\documentclass[aos,MSNbibl,seceqn,dvips]{arximspdf}
\usepackage{algorithmic,algorithm}
\usepackage{mathbh}
\usepackage{accents}

\doi{10.1214/13-AOS1101} 
\volume{41}
\issue{2}
\pubyear{2013}
\firstpage{693}
\lastpage{721}

\makeatletter

\floatname{algorithm}{Policy}

\renewcommand{\underline}{\underaccent{\bar}}

\newcommand{\textttt}[1]{{\fontsize{10pt}{11pt}\selectfont{\texttt{#1}}}}

\newcommand{\rrvert}{\vert}
\newcommand{\llvert}{\vert}

\newtheorem{theorem}{Theorem}[section]
\newtheorem{prop}{Proposition}[section]
\newtheorem{cor}{Corollary}[section]
\newtheorem{lem}{Lemma}[section]

\newcommand{\cA}{\mathcal{A}}
\newcommand{\cB}{\mathcal{B}}
\newcommand{\cC}{\mathcal{C}}
\newcommand{\cD}{\mathcal{D}}
\newcommand{\cG}{\mathcal{G}}
\newcommand{\cI}{\mathcal{I}}
\newcommand{\cJ}{\mathcal{J}}
\newcommand{\cL}{\mathcal{L}}
\newcommand{\cM}{\mathcal{M}}
\newcommand{\cP}{\mathcal{P}}
\newcommand{\cT}{\mathcal{T}}
\newcommand{\cX}{\mathcal{X}}
\newcommand{\blog}{\overline{\log}}

\newcommand{\R}{\mathbb R}
\newcommand{\p}{\mathbb P}
\newcommand{\E}{\mathbb E}
\newcommand{\1}{\mathbh1}
\newcommand{\N}{\mathbb N}

\renewcommand{\S}{\hat{\mathbb{I}}}

\newcommand{\pp}{{\mathsf p}}

\newcommand{\eps}{\varepsilon}

\newcommand{\argmin}{\mathop{\arg\min}}
\newcommand{\argmax}{\mathop{\arg\max}}

\newcommand{\ud}{\mathrm{d}}

\newcommand{\hpi}{\hat\pi}
\newcommand{\tpi}{\tilde\pi}
\newcommand{\bpi}{\bar\pi}
\newcommand{\sfR}{{\mathsf R}}
\newcommand{\wsfR}{\tilde{\mathsf R}}

\newcommand{\burst}{\mathsf{burst}}

\newcommand{\ucI}{\underline{\mathcal{I}}_B}
\newcommand{\ocI}{\bar{\mathcal{I}}_B}

\newcommand{\hdi}{\hat{\Delta}_{i}}
\newcommand{\hdj}{\hat{\Delta}_{j}}

\makeatother

\begin{document}
\begin{frontmatter}

\title{The multi-armed bandit problem with covariates}
\runtitle{Bandits with covariates}

\begin{aug}
\author[A]{\fnms{Vianney} \snm{Perchet}\thanksref{t1}\ead[label=e1]{vianney.perchet@normalesup.org}}
\and
\author[B]{\fnms{Philippe} \snm{Rigollet}\corref{}\thanksref{t2}\ead[label=e2]{rigollet@princeton.edu}}
\runauthor{V. Perchet and P. Rigollet}
\affiliation{Universit\'e Paris Diderot and Princeton University}
\address[A]{LPMA, UMR 7599\\
Universit\'e Paris Diderot\\
175, rue du Chevaleret\\
75013 Paris\\
France\\
\printead{e1}}
\address[B]{Department of Operations Research\\
\quad and Financial Engineering\\
Princeton University\\
Princeton, New Jersey 08544\\
USA\\
\printead{e2}} 
\end{aug}

\thankstext{t1}{Supported in part by the ANR Grant ANR-10-BLAN-0112.}
\thankstext{t2}{Supported in part by the NSF Grants DMS-09-06424,
DMS-10-53987.}

\received{\smonth{6} \syear{2012}}
\revised{\smonth{1} \syear{2013}}

%
\begin{abstract}
We consider a multi-armed bandit problem in a setting where each arm
produces a noisy reward realization which depends on an observable
random \textit{covariate}. As opposed to the traditional \textit{static}
multi-armed bandit problem, this setting allows for dynamically
changing rewards that better describe applications where side
information is available. We adopt a nonparametric model where the
expected rewards are smooth functions of the covariate and where the
hardness of the problem is captured by a \textit{margin} parameter. To
maximize the expected cumulative reward, we introduce a policy called
Adaptively Binned Successive Elimination (\textsc{abse}) that adaptively
decomposes the global problem into suitably ``localized'' static bandit
problems. This policy constructs an adaptive partition using a variant
of the Successive Elimination (\textsc{se}) policy. Our results include
sharper regret bounds for the \textsc{se} policy in a static bandit
problem and minimax optimal regret bounds for the \textsc{abse} policy in
the dynamic problem.
\end{abstract}

%

\end{frontmatter}

\section{Introduction}
\label{SECintro}

The seminal paper~\cite{Rob52} introduced an important
class of sequential optimization problems, otherwise known as
multi-armed bandits. These models have since been used
extensively in such fields as statistics, operations research,
engineering, computer science and economics. The traditional
multi-armed bandit problem can be described as follows. Consider
$K \ge2$ statistical populations (arms), where at each point in time it
is possible to sample from (pull) only one of them and receive a
random reward dictated by the properties of the sampled
population. The objective is to
devise a sampling policy that maximizes expected cumulative
rewards over a finite
time horizon. The difference between the performance of a given
sampling policy and that of an oracle, that repeatedly samples
from the population with the highest mean reward, is called the
\textit{regret}. Thus, one can re-phrase the objective as minimizing
the regret.

When the populations being sampled are homogeneous, that is, when the
sequential rewards are independent and identically distributed
(i.i.d.) in each arm, the family of upper-confidence-bound (UCB)
policies, introduced in~\cite{LaiRob85}, incur a regret of order
$\log n$, where $n$ is the length of the time horizon, and no
other ``good'' policy can (asymptotically) achieve a smaller
regret; see also~\cite{AueCesFis02}. The elegance of the theory
and sharp results developed in~\cite{LaiRob85} hinge to a large
extent on the assumption of homogenous populations and hence
identically distributed rewards. This, however, is clearly too
restrictive for many applications of interest. Often, the decision
maker observes further information and based on that, a more
\textit{customized} allocation can be made. In such settings, rewards may
still be assumed to be
independent, but no longer identically
distributed in each arm. A particular way to encode this is to
allow for an exogenous variable (a covariate) that affects the
rewards generated by each arm at each point in time when this arm
is pulled.

Such a formulation was first introduced in
\cite{Woo79} under parametric assumptions and in a somewhat
restricted setting; see~\cite{GolZee09,GolZee10} and \cite
{WanKulPoo05} for
very different recent approaches to the study of
such bandit problems, as well as references therein for further
links to antecedent literature. The first work to venture outside
the realm of parametric modeling assumptions appeared in
\cite{YanZhu02}. In particular, the mean response in each arm,
conditionally on the covariate value, was assumed to follow a general
functional form; hence one can view their setting as a \textit{nonparametric}
bandit problem. They propose a variant of the $\varepsilon$-greedy
policy (see, e.g.,~\cite{AueCesFis02}) and show that the average
regret tends to zero as the time horizon $n$ grows to infinity.
However, it is unclear whether this policy
satisfies a more refined notion of optimality, insofar as the magnitude
of the regret
is concerned, as is the case for
UCB-type policies in traditional bandit problems. Such questions were
partially addressed in~\cite{RigZee10} where near-optimal bounds on
the regret are proved in the case of a two-armed bandit problem under
only two
assumptions on the underlying functional form that governs the
arms' responses. The first is a mild smoothness condition, and the
second is a so-called \textit{margin condition} that involves a
\textit{margin parameter} which encodes the ``separation'' between the functions
that describe the arms' responses.

The purpose of the present paper is to extend the setup of \cite
{RigZee10} to the $K$-armed bandit problem with covariates when $K$ may
be large. This involves a customized definition of the margin
assumption. Moreover, the bounds proved in~\cite{RigZee10} suffered
two deficiencies. First, they hold only for a limited range of values
of the margin parameter and second, the upper bounds and the lower
bounds mismatch by a logarithmic factor. Improving upon these results
requires radically new ideas. To that end, we introduce three policies:
\begin{longlist}[(3)]
\item[(1)] Successive Elimination (\textsc{se}) is dedicated to the
static bandit case. It is the cornerstone of the other policies that
deal with covariates. During a first phase, this policy explores the
different arms, builds estimates and eliminates sequentially suboptimal
arms; when only one arm remains, it is pulled until the horizon is
reached. A variant of \textsc{se} was originally introduced in \cite
{EveManMan06}. However, it was not tuned to minimize the regret as
other measures of performance were investigated in this paper. We prove
new regret bounds for this policy that improve upon the canonical
papers~\cite{LaiRob85} and~\cite{AueCesFis02}.
\item[(2)] Binned Successive Elimination (\textsc{bse}) follows a simple
principle to solve the problem with covariates.
It consists of grouping similar covariates into bins and then looks
only at the average reward over each bin. These bins are viewed as
indexing ``local'' bandit problems,
solved by the aforementioned \textsc{se} policy. We prove optimal regret
bounds, polynomial in the horizon but only for a restricted class of
difficult problems. For the remaining class of easy problems, the
\textsc{bse} policy is suboptimal.
%
\item[(3)] Adaptively Binned Successive Elimination (\textsc{abse})
overcomes a severe limitation of the naive \textsc{bse}. Indeed, if the
problem is globally easy (this is characterized by the margin
condition), the \textsc{bse} policy employs a fixed and too fine
discretization of the covariate space. Instead, the \textsc{abse} policy
partitions the space of covariates in a fashion that adapts to the
local difficulty of the problem: cells are smaller when different arms
are hard to distinguish and bigger when one arm dominates the other.
This adaptive partitioning allows us to prove optimal regrets bounds
for the whole class of problems.
\end{longlist}
The optimal polynomial regret bounds that we prove are much larger than
the logarithmic bounds proved in the static case. Nevertheless, it is
important to keep in mind that they are valid for a much more flexible
model that incorporates covariates. In the particular case where $K=2$
and the problem is \textit{difficult}, these bounds improve upon the
results of~\cite{RigZee10} by removing a logarithmic factor that is
idiosyncratic to the \textit{exploration vs. exploitation} dilemma
encountered in bandit problems. Moreover, it follows immediately from
the previous minimax lower bounds of~\cite{AudTsy07} and \cite
{RigZee10}, that these bounds are optimal in a minimax sense and thus
cannot be further improved. It reveals an interesting and somewhat
surprising phenomenon: the price to pay for the partial information in
the bandit problem is dominated by the price to pay for nonparametric
estimation. Indeed the bound on the regret that we obtain in the bandit
setup for $K=2$ is of the same order as the best attainable bound in
the \textit{full information} case, where at each round, the operator
receives the reward from only one arm but observes the rewards of both
arms. An important example of the full information case is sequential
binary classification.


Our policies for the problem with covariates fall into the family of
``plug-in'' policies
as opposed ``minimum contrast'' policies; a detailed account of the
differences and similarities between these two setups in the full
information case can be found in~\cite{AudTsy07}. Minimum
contrast type policies have already received some attention in the
bandit literature with side information, aka \textit{contextual
bandits}, in the papers~\cite{LanZha08} and also
\cite{KakShaTew08}. A related problem online convex
optimization with side information was studied in
\cite{HazMeg07}, where the authors use a discretization technique
similar to the one employed in this paper. It is worth noting
that the cumulative regret in these papers is defined in a weaker
form compared to the traditional bandit literature, since the
cumulative reward of a proposed policy is compared to that of the
best policy in a certain restricted class of policies. Therefore,
bounds on the regret depend, among other things, on the complexity of
said class of
policies. Plug-in type policies have received attention in the
context of the continuum armed bandit problem, where as the name
suggests there are uncountably many arms. Notable entries in that
stream of work are~\cite{LuPalPal10} and~\cite{Sli11}, who
impose a smoothness condition both on the space of arms and the
space of covariates, obtaining optimal regret bounds up to
logarithmic terms.

\section{Improved regret bounds for the static problem}
\label{SECStatic}

In this section, it will be convenient for notational purposes, to
consider a multi-armed bandit problem with $K+1$ arms.

We revisit the Successive Elimination (\textsc{se}) policy introduced
in~\cite{EveManMan06} in the traditional setup of multi-armed bandit
problems. As opposed to the more popular UCB policy (see, e.g., \cite
{LaiRob85,AueCesFis02}), it allows us in the next section, to
construct an adaptive partition that is crucial to attain optimal rates
on the regret for the dynamic case with covariates. In this section, we
prove refined regret bounds for the \textsc{se} policy that exhibit a
better dependence on the expected rewards of the arms compared to the
bounds for UCB that were derived in~\cite{AueCesFis02}. Such an
improvement was recently attempted in~\cite{AueOrt10} and also
in~\cite{AudBub10} for modified UCB policies and we compare these
results to ours below.

Let us recall the traditional setup for the static multi-armed bandit
problem; see, for example,~\cite{AueCesFis02}.
Let $\cI=\{1,\ldots, K+1\}$ be a\vspace*{1pt} given set of $K+1 \ge2$ arms.
Successive pulls of arm $i \in\cI$ yield rewards $Y_1^{(i)}, Y_2^{(i)},
\ldots$ that are i.i.d. random variables in $[0,1]$ with expectation
given by\vspace*{1pt} $\E[Y_t^{(i)}] = f^{(i)}\in[0,1]$.
Assume without loss of generality that $f^{(1)} \le\cdots\le
f^{(K+1)}$ so that $K+1$ is one of the best arms. For simplicity, we
further assume that the best arm is \textit{unique} since for the
\textsc{se} policy, having multiple optimal arms only improves the regret
bound. In the analysis, it is convenient to denote this optimal arm by
$*:=K+1$ and to define the \textit{gaps} traditionally denoted by
$\Delta_1\ge\cdots\ge\Delta_*=0$, by $\Delta_i=f^{(*)}-f^{(i)}\ge0$.

A \textit{policy} $\pi=\{\pi_t\}$ is a sequence
of random variables $\pi_t \in\{1,\ldots, K+1\} $ indicating which
arm to pull at each time $t=1, \ldots, n$, and such that $\pi_t$
depends only on observations strictly anterior to $t$.

The performance of a policy
$\pi$ is measured by its
(\textit{cumulative}) \textit{regret} at time $n$ defined by
\[
R_n(\pi):= \sum_{t=1}^n
\bigl(f^{(*)}-f^{(\pi_t)} \bigr).
\]
Note that for a data-driven policy $\hpi$, this quantity is random
and, in the rest of the paper, we provide upper bounds on $\E R(\hpi
)$. Such bounds are referred to as \textit{regret bounds}.

We begin with a high-level description of the \textsc{se} policy
denoted by $\hpi$. It operates in rounds that are different from the
decision times $t=1,\ldots, n$. At the beginning of each round $\tau
$, a subset of the arms has been eliminated and only a subset $\cI
_\tau$ 
remains. During round $\tau$, each arm in $\cI_\tau$ is pulled
exactly once (\textsc{Exploration}). At the end of the round, for each
remaining arm in $\cI_\tau$, we decide whether to eliminate it using
a simple statistical hypothesis test: if we conclude that its mean is
significantly smaller than the mean of any remaining arm, then we
eliminate this arm and we keep it otherwise (\textsc{Elimination}). We
repeat this procedure until $n$ pulls have been made. The number of
rounds is random but obviously
smaller than~$n$.

The \textsc{se} policy, which is parameterized by two quantities $T \in
\N
$ and $\gamma>0$ and described in Policy~\ref{ALGse}, outputs an
infinite sequence of arms $\hpi_1, \hpi_2, \ldots$ without a
prescribed horizon. Of course, it can be truncated at any horizon $n$.
This description emphasizes the fact that the policy can be implemented
without perfect knowledge of the horizon $n$ and in particular, when
the horizon is a random variable with expected value $n$; 
nevertheless, in the static case, it is manifest from our result that,
when the horizon is known to be $n$, choosing $T=n$ is always the best
choice when possible and that other choices may lead to suboptimal results.

\begin{algorithm}[t]
\caption{Successive Elimination (\textsc{se})}
\label{ALGse}
\begin{algorithmic}
\REQUIRE Set of arms $\cI=\{1,\ldots, K\}$; parameters $T, \gamma$;
horizon $n$.
\ENSURE$(\hpi_1, \hat\tau_1, \S_1), (\hpi_2, \hat\tau_2, \S
_2), \ldots\in\cI\times\N\times\cP(\cI)$.
\STATE$\tau\leftarrow1$, $S \leftarrow\cI$, $t\leftarrow0$, $\bar
Y\leftarrow(0,\ldots, 0) \in[0,1]^K$
\LOOP

\STATE$\bar Y^{\max} \leftarrow\max\{\bar{Y}^{(i)}\dvtx i \in S\}$
\FOR{$i \in S$}
\IF{$\bar{Y}^{(i)} \ge\bar Y^{\max} - \gamma U(\tau, T)$}
\STATE$t \leftarrow t+1$
\STATE$\hpi_t \leftarrow i$ (observe $Y^{(i)}$) \hfill\textsc
{Exploration} 

\STATE$\S_t \leftarrow S$, $\hat\tau_t \leftarrow\tau$
\STATE$\bar Y^{(i)}\leftarrow\frac{1}{\tau}[( \tau-1) \bar
Y^{(i)}+{Y}^{(i)}]$
\ELSE
\STATE$S \leftarrow S \setminus\{i\}$. \hfill\textsc{Elimination}
\ENDIF
\ENDFOR
\STATE$\tau\leftarrow\tau+1$.
\ENDLOOP
\end{algorithmic}
\end{algorithm}

Note that after the exploration phase of each round $\tau=1, 2,\ldots,
$ each remaining arm $i \in\cI_\tau$ has been pulled exactly $\tau$
times, generating rewards $Y^{(i)}_1,\ldots, Y^{(i)}_\tau$. Denote by
$\bar{Y}^{(i)}(\tau)$ the average reward collected from arm $i \in
\cI_\tau$ at round $\tau$ that is defined by $\bar{Y}^{(i)}(\tau
)=(1/\tau) \sum_{t=1}^\tau Y^{(i)}_\tau$,
where here and throughout this paper, we use the convention $1/0=\infty
$. In the rest of the paper, $\log$ denotes the natural logarithm and
$\blog(x)=\log(x)\vee1$. For any positive integer $T$, define also
%
%
\begin{equation}
\label{EQdefU}
U(\tau,T)=2\sqrt{\frac{2\blog({T}/{\tau} )}{\tau
}},
\end{equation}
which is essentially a high probability upper bound on the magnitude of
deviations of $\bar{Y}^{(j)}(\tau)-\bar{Y}^{(i)}(\tau)$ from its
mean $f^{(j)}-f^{(i)}$.

The \textsc{se} policy for a $K$-armed bandit problem can be implemented
according to the pseudo-code of Policy~\ref{ALGse}. Note that, to
ease the presentation of Sections~\ref{SECbse} and~\ref{SECabse},
the \textsc{se} policy also returns at each time $t$, the number of rounds
$\hat\tau_t$ completed at time $t$ and a subset $\S_t \in\cP(\cI
)$ of arms that are active at time $t$, where $\cP(\cI)$ denotes the
power set of $\cI$.

The following theorem gives a first upper bound on the expected regret
of the \textsc{se} policy.

%
\begin{theorem}\label{THstatic}
Consider a $(K+1)$-armed bandit problem where horizon is a random
variable $N$ of expectation $n$ that is independent of the random
rewards. When implemented with parameters $T, \gamma\geq1$, the \textsc
{se} policy $\hpi$ exhibits an expected regret bounded, for any $\Delta
\ge0$, as
\[
\E\bigl[R_N(\hpi)\bigr]\leq392 \gamma^2 \biggl(1+
\frac{n}{T} \biggr)\frac
{K}{\Delta}\blog\biggl(\frac{T\Delta^2}{18\gamma^2} \biggr) +
n \Delta^-,
\]
%
where $\Delta^-$ is the largest $\Delta_j$ such that $\Delta
_j<\Delta$ if it exists, otherwise $\Delta^-=0$.
\end{theorem}
\begin{pf}
Assume without loss of generality that $\Delta_j >0$, for $j \ge1$
since arms $j$ such that $j=0$ do not contribute to the regret. Define
$\varepsilon_\tau=U(\tau,T)$. Moreover, for any $i$ in the set $\cI
_\tau$ of arms that remain active at the beginning of round $\tau$,
define $\hdi(\tau):=\bar Y^{(*)}(\tau) -\bar Y^{(i)}(\tau)$.
Recall that, at round $\tau$, if arms $i, * \in\cI_\tau$, then (i)
the optimal arm\vspace*{1pt} $*$ eliminates arm $i$ if $\hdi(\tau) \geq\gamma
\varepsilon_{\tau}$, and (ii) arm $i$ eliminates arm $*$ if $\hdi
(\tau) \leq- \gamma\varepsilon_{\tau}$.

Since $\hdi(\tau)$ estimates $\Delta_i$, the event in (i)
happens approximately, when $\gamma\varepsilon_\tau\simeq\Delta
_i$, so we introduce the
deterministic, but unknown, quantity\vadjust{\goodbreak} $\tau^*_i$ (and its approximation
$\tau_i=\lceil\tau^*_i \rceil$) defined as the solution of
%
\[
\Delta_i=\frac{3}{2}\gamma\varepsilon_{\tau^*_i}=3
\gamma\sqrt{\frac{2}{\tau^*_i}\blog\biggl(\frac{T}{\tau^*_i} \biggr)}
\qquad\mbox{so that } \tau_i \leq\tau_{i}^* +1 \leq
\frac{18\gamma^2
}{\Delta^2_{i}}\blog\biggl(\frac{T\Delta_{i}^2}{18\gamma^2} \biggr) +1.
\]
Note that $1\le\tau_1\le\cdots\le\tau_{K}$ as well as the bound
%
%
\begin{equation}
\label{EQboundtau} \tau_i\le\frac{19\gamma^2}{\Delta_i^2}\blog
\biggl(\frac{T\Delta
^2_i}{18\gamma^2} \biggr).
\end{equation}

We are going to decompose the regret accumulated by a suboptimal
arm $i$ into three quantities:
\begin{itemize}[--]
\item[--]the regret accumulated by pulling this arm at most until
round $\tau_i$: this regret is smaller than $\tau_i\Delta_i$;
\item[--]the regret accumulated by eliminating the optimal arm $*$
between round $\tau_{i-1}+1$ and $\tau_i$;
\item[--]the regret induced if arm $i$ is still present at round
$\tau_i$ (and in particular, if it has not been eliminated by the
optimal arm $*$).
\end{itemize}

We prove that the second and third events happen with small
probability, because of the choice of $\tau_i$. Formally, define
the following \textit{good} events:
\begin{eqnarray*}
\cA_i &=& \{ \mbox{the arm $*$ has not been eliminated
before round } \tau_i \};
\\
\cB_i &= & \bigl\{ \mbox{every arm } j \in\{1,\ldots,i\} \mbox
{ has been eliminated before round } \tau_j \bigr\}.
\end{eqnarray*}
Moreover, define $\cC_i=\cA_i\cap\cB_i$ and observe that $\cC_1
\supseteq\cC_2 \supseteq\cdots\supseteq C_K$. For any $i=1,\ldots,
K$, the contribution to the regret incurred after time $\tau_i$ on
$\cC_i$ is at most $N\Delta_{i+1}$ since each pull of arm $j \ge
i+1$ contributes to the regret by $\Delta_j \le\Delta_{i+1}$. We
decompose the underlying sample space denoted by $\cC_0$ into the
disjoint union $(\cC_0 \setminus\cC_{1})\cup\cdots\cup(\cC
_{K_0-1} \setminus\cC_{K_0})\cup\cC_{K_0}$ where $K_0 \in\{
1,\ldots,K\}$ is chosen later.
It implies the following decomposition of the expected regret:
%
%
\begin{equation}
\label{EQRegret1} \E R_N(\hat{\pi})\leq\sum
_{i=1}^{K_0}n\Delta_i\p(\cC
_{i-1}\setminus\cC_i )+\sum_{i=1}^{K_0}
\tau_i\Delta_i +n\Delta_{K_0+1}.
\end{equation}
Define by $A^c$ the complement of an event $A$. Note that the first
term on the right-hand side of the above inequality can be decomposed
as follows:
%
%
\begin{eqnarray}
\label{EQdecompCi} \sum_{i=1}^{K_0}n
\Delta_i\p(\cC_{i-1}\setminus\cC_i ) &=& n\sum
_{i=1}^{K_0}\Delta_i \p\bigl(
\cA_i^c \cap\cC_{i-1}\bigr)\nonumber\\[-8pt]\\[-8pt]
&&{}+n\sum
_{i=1}^{K_0}\Delta_i \p\bigl(
\cB_i^c\cap\cA_i\cap\cB_{i-1}
\bigr),\nonumber
\end{eqnarray}
where the right-hand side was obtained using the decomposition $\cC
_i^c=\cA_i^c \cup(\cB_i^c \cap\cA_i)$ and the fact that $\cA
_{i}\subseteq\cA_{i-1}$.

From Hoeffding's inequality, we have that for every $\tau\ge1$,
%
%
\begin{eqnarray}
\label{EqHoeff} \p\bigl(\hdi(\tau)< \gamma\varepsilon_\tau\bigr)&=&\p
\bigl(\hdi(\tau)-\Delta_i< \gamma\varepsilon_\tau-
\Delta_i \bigr)\nonumber\\[-8pt]\\[-8pt]
&\leq&\exp\biggl(-\frac{\tau(\Delta_i-\gamma\varepsilon
_\tau
)^2}{2} \biggr).\nonumber
\end{eqnarray}
On the event $\cB_i^c\cap\cA_i\cap\cB_{i-1}$, arm $*$ has not
eliminated arm $i$ at $\tau_i$. Therefore $\p(\cB_i^c\cap\cA_i\cap
\cB_{i-1}) \le\p(\hdi(\tau_i)<\gamma\varepsilon_{\tau_i})$.
Together with the above display with \mbox{$\tau= \tau_i$}, it yields
%
%
\begin{equation}
\label{EQboundBi} \p\bigl(\cB_i^c\cap
\cA_i\cap\cB_{i-1}\bigr)\leq\exp\biggl(-
\frac{\tau
_i\gamma^2\varepsilon_{\tau_i}^2}{8} \biggr)\leq\biggl(\frac{1}{e}
\wedge\frac{\tau_i}{T}
\biggr)^{\gamma^2}\leq\frac{\tau_i}{T},
\end{equation}
where we used the fact that $ \Delta_i \geq(3/2)\gamma\varepsilon
_{\tau_i}$.

It remains to bound the first term in the right-hand side of (\ref
{EQdecompCi}).
On the event $\cC_{i-1}$, the optimal arm $*$ has not been eliminated
before round $\tau_{i-1}$, but every suboptimal arm $j \leq i-1$ has.
So the probability that there exists an arm $j \ge i$ that eliminates
$*$ between $\tau_{i-1}$ and $\tau_i$ can be bounded as
\begin{eqnarray*}
\p\bigl(\cA_i^c\cap\cC_{i-1}\bigr) &\leq& \p
\bigl(\exists(j,s), i\le j \le K, \tau_{i-1}+1 \leq s \leq
\tau_i; \hdj(s) \leq-\gamma\varepsilon_s\bigr)
\\
&\leq& \sum_{j=i}^K\p\bigl(\exists s,
\tau_{i-1}+1 \leq s \leq\tau_i; \hdj(s) \leq- \gamma
\varepsilon_s\bigr)
\\
&=&\sum_{j=i}^K
\bigl[\Phi_j(\tau_i)-\Phi_j(
\tau_{i-1}) \bigr],
\end{eqnarray*}
where
$\Phi_j(\tau) = \p(\exists s \leq\tau; \hdj(s) \leq-\gamma
\varepsilon_s )$. Using Lemma~\ref{LEMpeeling}, we get $\Phi
_j(\tau) \leq4\tau/T$.
This bound implies that
\begin{eqnarray*}
&&
\sum_{i=1}^{K_0}\Delta_i \p
\bigl(\cA_i^c\cap\cC_{i-1}\bigr) \\
&&\qquad\leq\sum
_{i=1}^{K_0}\Delta_i \sum
_{j=i}^K \bigl[ \Phi_j(
\tau_i)-\Phi_j(\tau_{i-1}) \bigr]
\\
&&\qquad\le\sum_{j=1}^K
\sum_{i=1}^{j\wedge K_0-1} \Phi_j(
\tau_i) (\Delta_i -\Delta_{i+1} )+ \sum
_{j=1}^K \Phi_{j\wedge K_0}(\tau
_{j\wedge K_0})\Delta_{j\wedge K_0}
\\
&&\qquad\leq\frac{4}{T} \sum_{j=1}^K
\sum_{i=1}^{j\wedge K_0-1} \tau_i (
\Delta_i -\Delta_{i+1} )+\frac{4}{T}\sum
_{j=1}^K \tau_{j\wedge K_0}\Delta_{j\wedge K_0}.
\end{eqnarray*}
Using (\ref{EQboundtau}) and $\Delta_{i+1}\leq\Delta_i$, the
first sum can be bounded as
\begin{eqnarray*}
\sum_{j=1}^K \sum
_{i=1}^{j\wedge K_0-1} \tau_i (
\Delta_i -\Delta_{i+1} ) &\leq& 19\gamma^2\sum
_{j=1}^K \sum
_{i=1}^{j\wedge K_0-1} \blog\biggl(\frac{T\Delta^2_i}{18\gamma
^2} \biggr)
\frac{\Delta_i -\Delta_{i+1}}{\Delta_i^2}
\\
&\leq& 19\gamma^2 \sum_{j=1}^K
\int_{\Delta_{j\wedge K_0}}^{\Delta
_1} \blog\biggl(\frac{Tx^2}{18\gamma^2}
\biggr)\,\frac{\mathrm
{d}x}{x^2}
\\
&\leq& 19\gamma^2\sum_{j=1}^K
\frac{1}{\Delta_{j\wedge K_0}} \biggl[\blog\biggl(\frac{T\Delta
^2_{j\wedge K_0}}{18\gamma^2} \biggr)+2 \biggr].
\end{eqnarray*}
The previous two displays together with (\ref{EQboundtau}) yield
\[
\sum_{i=1}^{K_0}\Delta_i \p
\bigl(\cA_i^c\cap\cC_{i-1} \bigr) \leq
\frac{304\gamma^2}{T}\sum_{j=1}^K
\frac{1}{\Delta_{j\wedge
K_0}}\blog\biggl(\frac{T\Delta^2_{j\wedge K_0}}{18\gamma^2} \biggr).
\]
%
Putting together (\ref{EQRegret1}), (\ref{EQdecompCi}), (\ref
{EQboundBi}) and the above display
yield that the expected regret $\E R_N(\hpi)$of the \textsc{se} policy is
bounded above by
%
%
\begin{eqnarray}
\label{EQRegretgen}
&&
323\gamma^2 \biggl(1+ \frac{n}{T} \biggr)
\sum_{i=1}^{K_0} \frac{1}{\Delta_{i}}\blog
\biggl(\frac{n\Delta^2_{i}}{18\gamma^2} \biggr)\nonumber\\[-8pt]\\[-8pt]
&&\qquad{}+304\frac{\gamma^2n}{T}
\frac{K-K_0}{\Delta_{K_0}}\blog
\biggl(\frac{n\Delta^2_{K_0}}{18\gamma^2}
\biggr)+n\Delta_{K_0+1}.\nonumber
\end{eqnarray}
Fix $\Delta\ge0$ and let $K_0$ be such that $\Delta_{K_0+1}=\Delta
^-$. An easy study of the variations of the function
\[
x \mapsto\phi(x)=\frac{1}{x} \blog\biggl(\frac{nx^2}{18\gamma
^2} \biggr),\qquad
x > 0,
\]
reveals that $\phi(x) \le(2e^{-1/2}) \phi(x')$ for any $x \ge x'\ge0$.
%
%
%
Using this bound equation (\ref{EQRegretgen}) with $x'=\Delta_i,
i\le K_0$ and $x=\Delta$ completes the proof.
\end{pf}

The following corollary is obtained from a slight variations on the
proof of Theorem~\ref{THstatic}. It allows us to better compare our
results to the extant literature.

%
\begin{cor}\label{CorStatic}
Under the setup of Theorem~\ref{THstatic}, the \textsc{se} policy $\hpi
$ run with parameter $T=n$ and $\gamma=1$ satisfies for any $K_0 \leq K$,
%
%
\begin{equation}
\label{EQCORSTATIC1}\quad \E R_N(\hpi)\leq646\sum
_{i=1}^{K_0} \frac{\blog(n\Delta
^2_{i} )}{\Delta_{i}}+304\frac{K-K_0}{\Delta_{K_0}}
\blog\bigl(n\Delta^2_{K_0} \bigr)+n\Delta_{K_0+1}.
\end{equation}
In particular,
%
%
\begin{equation}
\label{EQTHstatic} \E R_N(\hpi)\leq\min\Biggl\{ 646\sum
_{i=1}^{K} \frac{\blog
(n\Delta^2_{i} )}{\Delta_{i}}, 166
\sqrt{nK\log(K)} \Biggr\}.
\end{equation}
\end{cor}
\begin{pf} Note that (\ref{EQCORSTATIC1}) follows from (\ref
{EQRegretgen}). To prove (\ref{EQTHstatic}), take $K_0=K$ in (\ref
{EQCORSTATIC1}) and $\Delta=28\sqrt{K\log(784K/18)/n}$
in Theorem~\ref{THstatic}, respectively.
\end{pf}
%
This corollary is actually closer to the result of~\cite{AueOrt10}.
The additional second term in our bound comes from the fact that we had
to take into account the probability that an optimal arm $*$ can be
eliminated by any arm, not just by some \textit{suboptimal arm} with
index lower than $K_0$; see~\cite{AueOrt10}, page 8. It is unclear why
it is enough to look at the elimination by those arms, since if $*$ is
eliminated---no matter the arm that eliminated it---the Hoeffding
bound (\ref{EqHoeff}) no longer holds.

The right-hand side of (\ref{EQTHstatic}) is the minimum of two
terms. The first term is distribution-dependent and shows that the
\textsc{se} policy adapts to the unknown distribution of the rewards.
It is
very much in the spirit of the original bound of~\cite{LaiRob85} and
of the more recent finite sample result of~\cite{AueCesFis02}. Our
bound for the \textsc{se} policy is smaller than the aforementioned bounds
for the UCB policy by a logarithmic factor. Reference~\cite{LaiRob85} did not
provide the first bounds on the expected regret. Indeed,~\cite{Vog60}
and~\cite{Bat81} had previously derived what is often called \textit
{gap-free} bound as they hold uniformly over the $\Delta_i$'s. The
second term in our bound is such a gap-free bound. It is of secondary
interest in this paper and arise as a byproduct of refined distribution
dependent bound. Nevertheless, it allows us to recover near optimal
bounds of the same order as~\cite{JudNazTsy08}. They depart from
optimal rates by a factor $\sqrt{\log K}$ as proved in~\cite
{AudBub10}. Actually, the result of~\cite{AudBub10} is much stronger
than our gap-free bound since it holds for any sequence of bounded
rewards, not necessarily drawn independently. 

None of the distribution-dependent bounds in Corollary \ref
{CorStatic} or the one provided in~\cite{AudBub10} is stronger than
the other. The superiority of one over the other depends on the set $\{
\Delta_1,\ldots, \Delta_K\}$: in some cases (e.g., if all
suboptimal arms have the same expectation) the latter is the best while
in other cases (if the $\Delta_i$ are spread) our bounds are better.




\section{Bandit with covariates}
\label{secdesc}

This section is dedicated to a detailed description of the
nonparametric bandit with covariates.

\subsection{Machine and game}\label{subbanditcov}
A \textit{$K$-armed bandit machine with covariates} (with $K$ an integer
greater than 2) is
characterized by a sequence
\[
\bigl(X_t, Y_t^{(1)}, \ldots,
Y_t^{(K)}\bigr),\qquad t=1,2,\ldots,
\]
of independent random vectors, where $ (X_t )_{t\ge1}$,
is a sequence of i.i.d. covariates in $\cX=[0,1]^d$ with probability
distribution $P_X$, and $Y_t^{(i)}$ denotes
the random reward yielded by arm $i$ at time $t$. Throughout the paper,
we assume that $P_X$ has a density, with respect to the Lebesgue
measure, bounded above and below by some $\bar{c}>0$ and
$\underline{c}>0$, respectively. We denote by $E_X$ the expectation
with respect to $P_X$. We assume that,
for each $i\in\cI=\{1,\ldots, K\}$, 
rewards $Y_t^{(i)}, t=1, \ldots, n$, are
random variables in $[0,1]$ with conditional expectation given by
\[
\E\bigl[Y_t^{(i)}|X_t\bigr] =
f^{(i)}(X_t),\qquad i=1,\ldots, K, t=1,2, \ldots,
\]
where $f^{(i)}, i=1,\ldots, K$, are unknown functions such that $0
\le
f^{(i)}(x) \le1$, for any $i=1, \ldots, K, x \in\cX$. A natural
example is
where $Y_t^{(i)}$ takes values in $\{0,1\}$ so that the
conditional distribution of $Y_t^{(i)}$ given $X_t$ is Bernoulli with
parameter $f^{(i)}(X_t)$.

The \textit{game} takes place sequentially on this machine, pulling
one of the arms at each time $t=1,\ldots, n$. A
\textit{policy} $\pi=\{\pi_t\}$ is a sequence
of random functions $\pi_t\dvtx\cX\to\{1,\ldots, K\} $ indicating to the
operator which arm to pull at each time~$t$, and such that $\pi_t$
depends only on observations strictly anterior to $t$. The
\textit{oracle policy}~$\pi^\star$, refers to the strategy that would
be run by an omniscient operator with complete knowledge of the
functions $f^{(i)}, i=1, \ldots, K$. Given side information~$X_t$, the oracle
policy $\pi^{\star}$ prescribes to pull any arm with the largest
expected reward,
that is,
\[
\pi^\star(X_t) \in\argmax_{i=1,\ldots,K}
f^{(i)}(X_t)
\]
with ties broken arbitrarily. Note that the function $f^{(\pi^\star
(x))}(x)$ is equal to the pointwise maximum of the functions $f^{(i)},
i=1,\ldots, K$, defined by
\[
f^{\star}(x)=\max\bigl\{ f^{(i)}(x); i =1,\ldots,K \bigr\}.
\]
The oracle rule is used to benchmark any proposed policy
$\pi$ and to measure the performance of the latter via its
(\textit{cumulative}) \textit{regret} at time $n$ defined by
\[
R_n(\pi):=\E\sum_{t=1}^n
\bigl(Y^{(\pi^\star(X_t))}_t-Y^{(\pi
_t(X_t))}_t \bigr)= \sum
_{t=1}^n E_X
\bigl(f^\star(X)-f^{(\pi
_t(X))}(X) \bigr).
\]

Without further assumptions on the machine, the game can be
arbitrarily difficult and, as a result, expected regret can be
arbitrarily close to $n$. In the following
subsection, we describe natural regularity conditions under which
it is possible to achieve sublinear growth rate of the expected regret,
and characterize policies that perform in
a near-optimal manner.

\subsection{Smoothness and margin conditions}

As usual in nonparametric estimation we first impose some
regularity on the functions $f^{(i)}, i=1,\ldots,K$. Here and in what
follows we use $\|\cdot\|$ to denote the Euclidean norm on $\R
^d$.\vadjust{\goodbreak}

\textsc{Smoothness condition.} We say\vspace*{1pt} that the machine
satisfies the smoothness condition with parameters $(\beta, L)$ if
$f^{(i)}$ is $(\beta, L)$-H\"older, that is, if
\[
\bigl|f^{(i)}(x)-f^{(i)}\bigl(x'\bigr)\bigr|\le L
\bigl\|x-x'\bigr\|^\beta\qquad \forall x, x' \in\cX, i=1,
\ldots,K,
\]
for some $\beta\in(0,1]$ and $L>0$.\vspace*{9pt}

Now denote the second pointwise maximum of the functions $f^{(i)},
i=1,\ldots, K$, by $f^{\sharp}$; formally for every $x \in\cX$ such that
$\min_i f^{(i)}(x) \neq\max_i f^{(i)}(x)$ it is defined by
\[
f^{\sharp}(x)= \max_i \bigl\{ f^{(i)}(x);
f^{(i)}(x) < f^\star(x) \bigr\}
\]
and by $f^{\sharp}(x)=f^{\star}(x)=f^{(1)}(x)$ otherwise. Notice that
a direct consequence of the smoothness condition is that the function
$f^{\star}$ is $(\beta, L)$-H\"older; however, $f^\sharp$~might not
even be continuous.

The behavior of the function
$\Delta:=f^\star-f^\sharp$ critically controls the complexity of the
problem and the H\"older regularity gives a local upper bound on this quantity.
The second condition gives a lower bound on this function though in a weaker
global sense. It is closely related to the margin condition
employed in classification~\cite{Tsy04,MamTsy99}, which drives the
terminology employed here. It was originally imported to the bandit setup
by~\cite{GolZee09}.\vspace*{9pt}

\textsc{Margin condition.} We say that the machine
satisfies the margin condition with parameter $\alpha>0$ if there
exists $\delta_0\in(0,1)$, $C_0>0$ such that
\[
P_X \bigl[ 0< f^{\star}(X)-f^{\sharp}(X) \le\delta
\bigr]\le C_0 \delta^\alpha\qquad \forall\delta\in[0,
\delta_0].
\]
%

If the marginal $P_X$
has a density bounded above and below, the margin condition contains only
information about the behavior of the function $\Delta$ and not
the marginal $P_X$ itself. This is in contrast with~\cite{GolZee09}
where the margin assumption is used precisely to control the behavior
of the marginal $P_X$ while that of the reward functions is fixed. A
large value of the parameter $\alpha$
means that the function $\Delta$ either takes value 0 or is
bounded away from 0, except over a set of small
$P_X$-probability. Conversely, for values of $\alpha$ close to 0,
the margin condition is essentially void, and the two functions can
be arbitrary close, making it difficult to distinguish them. This
reflects in the bounds on the expected regret 
derived in
the subsequent section.

Intuitively, the smoothness condition and the margin condition
work in opposite directions. Indeed, the former ensures that the
function $\Delta$ does not ``depart from zero'' too fast whereas
the latter warrants the opposite. The following proposition
quantifies the extent of this conflict.

%
\begin{prop}
\label{propalpha-beta} Under the smoothness condition with
parameters\break $(\beta, L)$, and the margin condition with parameter
$\alpha$, the following hold:
\begin{itemize}[--]
\item[--]if $\alpha\beta>d$, then a given arm is either always or
never optimal; in the latter case, it is bounded away from $f^\star$
and one can take $\alpha=\infty$;
\item[--]if $\alpha\beta\le d$, then there exist machines with
nontrivial oracle policies.
\end{itemize}
\end{prop}
\begin{pf} This proposition is a straightforward consequences of,
respectively, the first two points of Proposition 3.4 in~\cite{AudTsy05}.

For completeness, we provide an example with $d=1$,
$\cX=[0,1]$, $f^{(2)}=\cdots=f^{(K)}\equiv0$ and $f^{(1)}(x)=L
\operatorname{sign}(x-0.5)|x-0.5|^{1/\alpha}$. Notice that $f^{(1)}$ is
$(\beta,
L)$-H\"older if and only if
$\alpha\beta\le1$. Any oracle policy is nontrivial, and, for
example, one can define $\pi^{\star}(x)=2$ if $x\le0.5$ and $\pi
^{\star}(x)=1$ if
$x>0.5$. Moreover, it can be easily shown that the machine satisfies
the margin condition with parameter $\alpha$ and with $\delta_0=C_0=1$.
\end{pf}

We denote by $\cM_{\cX}^K(\alpha,\beta,L)$ the class of $K$-armed
bandit problems with covariates in $\cX=[0,1]^d$ with a machine
satisfying the margin condition with parameter $\alpha>0$, the
smoothness condition with parameters $(\beta,L)$ and such that $P_X$
has a density, with respect to the Lebesgue measure, bounded above and
below by some $\bar{c}>0$ and $\underline{c}>0$, respectively.

\subsection{Binning of the covariate space}
\label{subbin}

To design a policy that solves the bandit problem with covariates
described above, one has to inevitably find an estimate of the
functions $f^{(i)}, i=1,\ldots,K$, at the current point $X_t$. There
exists a
wide variety of nonparametric regression estimators ranging from
local polynomials to wavelet estimators. Both of the policies
introduced below are based on estimators of $f^{(i)}, i=1,\ldots,K$, that
are $P_X$ almost surely piecewise constant over a particular collection
of subsets, called \textit{bins} of the covariate space $\cX$.

We define a partition of $\cX$ in a measure theoretic sense as a
collection of measurable sets, hereafter called \textit{bins},\vspace*{1pt} $B_1,
B_2, \ldots$ such that $P_X(B_j)>0$,
$\bigcup_{j\ge1} B_j=\cX$ and $B_j \cap B_k =\varnothing$, $ j,k\ge
1$, up to sets of null $P_X$ probability.
For any $i \in\{\star, 1,\ldots, K\}$ and any bin $B$, define
%
%
\begin{equation}
\label{EQaverage} \bar f^{(i)}_B=\E\bigl[f^{(i)}(X_t)
| X_t \in B\bigr]=\frac{1}{P_X(B)}\int_{B}f^{(i)}
(x)\,\ud P_X(x).
\end{equation}
To define and analyze our policies, it is convenient to reindex the
random vectors $(X_t, Y_t^{(1)}, \ldots, Y_t^{(K)})_{t\ge1}$ as
follows. Given a bin $B$, let $t_{B}(s)$ denote the $s$th time at which
the sequence $(X_t)_{t\ge1}$ is in $B$ and observe that it is a
stopping time. It is a standard exercise to show that, for any bin $B$
and any arm $i$, the random variables $Y_{t_{B}(s)}^{(i)}$, $s \ge1$ are
i.i.d. with expectation given by $\bar f_B^{(i)}\in[0,1]$. As a result,
the random variables $Y_{B,1}^{(i)}, Y_{B,2}^{(i)}, \ldots$ obtained by
successive pulls of arm $i$ when $X_t \in B$ form an i.i.d. sequence in
$[0,1]$ with expectation given by $\bar f_B^{(i)}\in[0,1]$.
Therefore, if we restrict our attention to observations in a given bin
$B$, we are in the same setup as the static bandit problem studied in
the previous section. This observation leads to the notion of \textit
{policy on $B$}. More precisely, fix a subset $B \subset\cX$, an
integer $t_0\ge1$ and recall that $\{ t_B(s)\dvtx s\ge1, t_B(s) \ge
t_0\}$ is the set of chronological times $t$ posterior to $t_0$ at
which $X_t \in B$. Fix $\cI' \subset\cI$ and consider the static
bandit problem with arms $\cI'$ defined in Section~\ref{SECStatic}
where successive pulls of arm $i \in\cI'$, at times posterior to
$t_0$, yield rewards $Y^{(i)}_{B,1}, Y^{(i)}_{B,2},\ldots,$ that are
i.i.d. in $[0,1]$ with mean $\bar f_B^{(i)}\in[0,1]$. The \textsc{se}
policy with parameters $T,\gamma$ on this static problem is called
\textit{\textsc{se} policy on $B$ initialized at time $t_0$ with initial
set of arms $\cI'$ and parameters $T, \gamma$}.


%

\section{Binned Successive Elimination}
\label{SECbse}

We first outline a naive policy to operate the
bandit machine described in Section~\ref{secdesc}. It consists of
fixing a partition of $\cX$ and for each set $B$ in this partition, to
run the \textsc{se} policy on $B$ initialized at time $t_0=1$ with initial
set of arms $\cI$ and parameters $T,\gamma$ to be defined below.

The \textit{Binned Successive Elimination} (\textsc{bse}) policy $\bpi$
relies on a specific partition of $\cX$. Let $\cB_M:=\{B_1,\ldots,
B_{M^d}\}$ be the regular partition of $\cX= [0,1]^d$: the collection of
hypercubes defined for ${\mathsf k}=(k_1,\ldots, k_d)\in\{1,\ldots, M\}
^d$ by
\[
B_{\mathsf k}= \biggl\{x\in\cX\dvtx \frac{k_\ell-1}{M}\le x_\ell\le
\frac{k_\ell}{M}, \ell=1,\ldots, d \biggr\}.
\]
In this paper, all sets are defined up to sets of null Lebesgue
measure. As mentioned in Section~\ref{subbin}, the problem can be
decomposed into $M^d$ independent static bandit problems, one for each
$B \in\cB_M$.

%
Denote by $\hpi_B$ the \textsc{se} policy on bin $B$
with initial set of arms $\cI$ and parameters $T=nM^{-d}, \gamma=1$.
For any $x \in\cX$, let $\cB(x) \in\cB_M$ denote the bin such that
$x \in\cB(x)$. Moreover, for any time $t \ge1$, define
%
%
\begin{equation}
\label{EQdefNBj} N_B(t)=\sum_{l=1}^t
\1(X_l \in B)
\end{equation}
to be the number of times before $t$ when the covariate fell in bin
$B$. The \textsc{bse} policy $\bar\pi$ is a sequence of functions $\bar
\pi_t\dvtx\cX\to\cI$ defined by $\bar\pi_t(x)=\hpi_{B,N_{B}(t)}$,
where $B=\cB(x)$. It can be implemented according to the pseudo-code
of Policy~\ref{ALGbse}.

\begin{algorithm}[t]
\caption{Binned Successive Elimination (\textsc{bse})}
\label{ALGbse}
\begin{algorithmic}
\REQUIRE Set of arms $\cI=\{1,\ldots, K\}$. Parameters $n, M$.
\ENSURE$\bpi_1, \ldots,\bpi_n \in\cI$.
\STATE$\cB\leftarrow\cB_M$
\FOR{$B \in\cB_M$}
\STATE Initialize a \textsc{se} policy $\hpi_B$ with parameters
$T=nM^{-d}, \gamma=1$.
\STATE$N_B\leftarrow0$.
\ENDFOR
\smallskip
\FOR{$t=1,\ldots, n$}
\STATE$B \leftarrow\cB(X_t)$.
\STATE$N_{B} \leftarrow N_B + 1$.
\STATE$\bpi_t \leftarrow\hpi_{B,N_B}$ (observe $Y^{(\bpi_t)}_t$).
\ENDFOR
\end{algorithmic}
\end{algorithm}

The following theorem gives an upper bound on the expected regret of
the \textsc{bse} policy in the case where the problem is difficult, that
is, when the margin parameter $\alpha$ satisfies $0<\alpha<1$.

%
\begin{theorem}
\label{THmain} Fix $\beta\in(0,1]$, $L>0$ and $\alpha\in
(0,1)$ and consider a problem in $\cM^K_{\cX}(\alpha,\beta,L)$.
Then the \textsc{bse} policy $\bpi$
with $M=\lfloor(\frac{n}{K\log(K)} )^{1/(2\beta
+d)}\rfloor$ has an expected regret at time $n$ bounded as follows:
\[
\E R_n(\bpi)\le C n \biggl(\frac{K\log K}{n} \biggr)^{{\beta
(\alpha+1)}/({2\beta+d})},
\]
where $C>0$ is a positive constant that does not depend on $K$.
\end{theorem}
The case $K=2$ was studied in~\cite{RigZee10} using a similar policy
called UCBogram. Unlike in~\cite{RigZee10} where suboptimal bounds for
the UCB policy are used, we use here the sharper regret bounds of
Theorem~\ref{THstatic} and the \textsc{se} policy as a running horse for
our policy, thus leading to a better bound than~\cite{RigZee10}.
Optimality for the two-armed case is discussed after Theorem~\ref{THmainadap}.

\begin{pf*}{Proof of Theorem~\ref{THmain}}
We assume that $\cB_M=\{B_1,\ldots, B_{M^d}\}$ where
the indexing will be made clearer later in the proof. Moreover, to keep
track of positive constants, we number them $c_1, c_2, \ldots\,$. For
any real valued function $f$ on $\cX$ and any measurable $A \subseteq
\cX$, we use the notation $P_X(f \in A)=P_X(f(X) \in A)$. Moreover,
for any $i \in\{\star, 1,\ldots, K\}$, we use the notation $\bar
f_j^{(i)}=\bar f_{B_j}^{(i)}$.

Define $c_1=2Ld^{\beta/2} +1$, and let $n_0 \ge2$ be the largest
integer such that
$n_0^{\beta/(2\beta+d)} \le2c_1/\delta_0$,
where $\delta_0$ is the constant appearing in the margin condition.
If $n \le n_0$, we have $R_{n}(\bar\pi) \le n_0$ so that the result
of the theorem holds when $C$ is chosen large enough, depending on the
constant $n_0$.
In the rest of the proof, we assume that $n > n_0$ so that
$c_1M^{-\beta} < \delta_0$.

Recall that the \textsc{bse}
policy $\bar\pi$ is a collection of functions $\bar\pi_{t}(x)=\hat
\pi_{\cB(x),N_{\cB(x)}(t)}$ that are constant on each
$B_j$. Therefore, the regret of $\bar\pi$ can be decomposed as
$R_n(\bar\pi)=
\sum_{j=1}^{M^d}\sfR_j(\bar\pi)$, where
%
\[
{\sfR_j(\bar\pi)=\sum_{t=1}^n
\bigl(f^{\star}(X_t)-f^{(\hat\pi
_{B,N_B(t)})}(X_t) \bigr)
\1(X_t \in B_j).}
\]
Conditioning on the event $\{X_t \in B_j\}$, it follows from (\ref
{EQaverage}) that
\[
\E\sfR_j(\bar\pi)=\E\Biggl[\sum_{t=1}^n
\bigl(\bar f_j^\star-\bar f_j^{(\bar\pi_t)}
\bigr)\1(X_t \in B_j) \Biggr]=\E\Biggl[\sum
_{s=1}^{N_{j}(n)} \bigl(\bar f_j^\star
-\bar f{}^{(\hat{\pi}_{B_j,s})}_j \bigr) \Biggr],
\]
where $N_{j}(n)=N_{B_j}(t)$ is defined in (\ref{EQdefNBj}); it
satisfies, by assumption, $\underline{c}nM^{-d} \leq\E
[N_{j}(n)]\leq\bar{c}nM^{-d}$.

Let $\wsfR_j(\bar\pi)=\sum_{s=1}^{N_{j}(n)} f^*_j-\bar f_j^{(\hat
{\pi}_{B_j,s})}$ be the regret associated to a static bandit problem
with arm $i$ yielding reward $\bar f^{(i)}_j$ and where $ f^*_j=\max
_i\bar f^{(i)}_j \le\bar f^\star_j$ is the largest average reward. It
follows from the smoothness condition that $\bar f^\star_j\le
f^*_j+c_1M^{-\beta}$ so that
%
%
\begin{equation}
\label{EQwsfR}\quad \E\sfR_j(\bar\pi)\le\E\wsfR_j(\bar
\pi)+\bar c nM^{-d}\bigl(\bar f^\star_j-
f^*_j\bigr)\le\E\wsfR_j(\bar\pi) + c_1
\bar c nM^{-\beta
-d}.
\end{equation}

Consider \textit{well-behaved} bins on which the expected reward
functions are well separated. These are bins $B_j$ with indices in $\cJ
$ defined by
\[
\cJ:=\bigl\{j\in\bigl\{1,\ldots,M^d\bigr\} \mbox{ s.t. } \exists x
\in B_j, f^{\star}(x)-f^{\sharp}(x)>
c_1M^{-\beta} \bigr\}.
\]
A bin $B$ that is not well behaved is called \textit{strongly ill
behaved} if there is some $x \in B$ such that $f^\star(x)=f^\sharp(x)
= f^{(i)}(x)$, for all $i \in\cI$, and \textit{weakly ill behaved}
otherwise. Strongly and weakly ill behaved bins have indices in
\[
\cJ_s^c:= \bigl\{j\in\bigl\{1,\ldots,M^d
\bigr\} \mbox{ s.t. } \exists x \in B_j, f^\star(x)=f^\sharp(x)
\bigr\}
\]
and
\[
\cJ^c_w:=\bigl\{j\in\bigl\{1,\ldots,M^d
\bigr\} \mbox{ s.t. } \forall x \in B_j, 0<f^{\star}(x)-f^{\sharp}(x)
\le c_1M^{-\beta} \bigr\},
\]
%
respectively. Note that for any $i \in\cI$, the function $f^\star
-f^{(i)}$ is $(\beta, 2L)$-H\"older. Thus for any $j \in\cJ_s^c$ and
any $i = 1,\ldots, K$, we have $f^\star(x)-f^{(i)}(x) \le
c_1M^{-\beta
}$ for all $x \in B_j$ so that the inclusion $\cJ^c_s \subset\{
1,\ldots,M^d\}\setminus\cJ$ indeed holds.\vspace*{9pt}

\textit{First part}: \textit{Strongly ill behaved bins in $\cJ^c_s$.\quad}
Recall that for any $j \in\cJ^c_s$, any arm $i\in\cI$, and any $x
\in B_j$, $f^{\star}(x) - f^{(i)}(x) \leq c_1M^{-\beta}$. Therefore,
%
%
\begin{eqnarray}\label{eqstronglyill}
\sum_{j\in\cJ_s^c} \E\sfR_j(\bpi)
&\le&
c_1 n M^{-\beta}P_X \bigl\{0<f^\star(X)-f^\sharp(X)
\leq c_1M^{-\beta} \bigr\}
\nonumber\\[-8pt]\\[-8pt]
&\leq& c_1^{1+\alpha} n M^{-\beta(1+\alpha)},\nonumber
\end{eqnarray}
where we used the fact that the set $\{x \in\cX\dvtx f^\star
(x)=f^\sharp(x)\}$ does not contribute to the regret.
\vspace*{9pt}

\textit{Second part}: \textit{Weakly ill behaved bins in $\cJ^c_w$.\quad}
The numbers\vspace*{1pt} of weakly ill behaved bins can be bounded using $f^\star
(x) - f^{\sharp}(x)>0$ on such a bin; indeed, the margin condition
implies that
\[
\sum_{j \in\cJ_w^c} \frac{\underline c}{M^d} \leq
P_X \bigl\{0< f^{\star
}(X)-f^{\sharp}(X)\leq
c_1M^{-\beta} \bigr\} \leq{c_1^\alpha}
M^{-\beta\alpha}.
\]
It yields $|\cJ_w^c| \le\frac{c_1^\alpha}{\underline c}M^{d-\beta
\alpha}$. Moreover, we bound the expected regret on weakly ill behaved
bins using Theorem~\ref{THstatic} with specific values
\[
\Delta^-<\Delta:=\sqrt{K\log(K)M^d/n},\qquad \gamma=1 \quad\mbox{and}\quad T=nM^{-d}.
\]
Together with (\ref{EQwsfR}), it yields
%
%
\begin{equation}
\label{eqweaklyill} \sum_{j\in\cJ_w^c} \E
\sfR_j(\bpi) \le
c_2 \bigl[
\sqrt{K\log(K)} M^{{d}/{2}-\beta\alpha} \sqrt{n}+nM^{-\beta
(1+\alpha)} \bigr].
\end{equation}

\textit{Third part}: \textit{Well-behaved bins in $\cJ$.\quad}
This part is decomposed into two steps. In the first step, we bound the
expected regret in a given bin $B_j, j \in\cJ$; in the second step we
use the margin condition to control the sum of all these expected regrets.

Step 1. Fix $j \in\cJ$ and recall that there exists $x_j \in B_j$
such that $f^{\star}(x_j)-f^{\sharp}(x_j)>c_1M^{-\beta}$. Define
$\cI_j^\star=\{i\in\cI\dvtx f^{(i)}(x_j)=f^\star(x_j)\}$ and $\cI
_j^0=\cI\setminus\cI_j^\star=\{i \in\cI\dvtx f^{\star}(x_j)-f^{(i)}
(x_j)>c_1M^{-\beta}\}$. We call $\cI_j^\star$ the set of (almost)
\textit{optimal} arms over $B_j$ and $\cI_j^0$ the set of \textit
{suboptimal} arms over $B_j$. Note that $\cI_j^0 \neq\varnothing$ for
any $j \in\cJ$.

The smoothness condition implies that for any $i \in\cI_j^0, x \in B_j$,
%
%
\begin{equation}
\label{EQInaught} f^\star(x) -f^{(i)}(x) >
c_1M^{-\beta} - 2L\|x-x_j\|^\beta\ge
M^{-\beta}.
\end{equation}
Therefore, $f^{\star}-f^{\sharp}>0$ on $B_j$. Moreover, for any arm
$i \in\cI_j^\star$ that is not the best arm at some $x \neq x_j$,
then necessarily $0<f^\star(x) - f^\sharp(x)\leq f^\star(x)-f^{(i)}(x)
\leq c_1M^{-\beta}$.
So for any $x \in B_j$ and any $i \in\cI_j^\star$, it holds that
either $f^\star(x)= f^{(i)}(x)$ or $f^\star(x) -f^{(i)}(x)\le
c_1M^{-\beta}$. It yields
%
%
\begin{equation}
\label{EQIstar} f^\star(x) -f^{(i)}(x)\le
c_1M^{-\beta}\1 \bigl\{ 0 < f^\star(x)
-f^\sharp(x) \le c_1M^{-\beta} \bigr\}.
\end{equation}
Thus, for any optimal arm $i \in\cI_j^\star$, the reward functions
averaged over $B_j$ satisfy $\bar f^\star_j -\bar f^{(i)}_j\le
c_1M^{-\beta}q_j$, where
\[
q_j:=P_X \bigl\{ 0 < f^\star-f^\sharp
\le c_1M^{-\beta} | X \in B_j \bigr\}.
\]
Together with (\ref{EQwsfR}), it yields $\E\wsfR_j(\bpi)\le\E
\sfR_j(\bpi)+\bar c c_1nM^{-d-\beta}q_j$.
For any suboptimal arms $i \in\cI_j^0$, (\ref{EQInaught}) implies
that $\underline{\Delta}_j^{(i)}:=\bar{f}^{\star}_j-\bar{f}^{(i)}_j
> M^{-\beta}$.

Assume now without loss of generality that the average gaps $\underline
{\Delta}_j^{(i)}$ are ordered in such a way that $\underline{\Delta
}_j^{(1)} \ge\cdots\ge\underline{\Delta}_j^{(K)}$. Define
\[
K_0:=\argmin_{i \in\cI_j^0} \underline{\Delta}_j^{(i)}
\quad\mbox{and}\quad \underline{\Delta}_j:=\underline{\Delta}_j^{(K_0)}
\]
and observe that if $i \in\cJ$ is such that $\underline{\Delta
}_j^{(i)}<\underline{\Delta}_j$, then $i \in\cI_j^\star$. Therefore,
it follows from (\ref{EQIstar}) that $\underline{\Delta}_j^{(i)}\le
c_1M^{-\beta}q_j$ for such $i$.
Applying Theorem~\ref{THstatic} with $\underline{\Delta}_j$ as
above and $\gamma=1$, we find that there exists a constant $c_3>0$
such that, for any $j\in\cJ$,
\[
\E\wsfR_j(\bpi)\le392(1+\bar c)\frac{K}{ \underline{\Delta
}_j}\blog
\bigl(nM^{-d}\underline{\Delta}_j^2 \bigr)+
\bar c c_1nM^{-d-\beta}q_j.
\]
Hence,
%
%
\begin{equation}
\label{EQregJ} \E\sfR_j(\bpi)\le c_3 \biggl(
\frac{K}{ \underline{\Delta}_j}\blog\bigl(nM^{-d}\underline{\Delta}_j^2
\bigr)+nM^{-d-\beta}q_j \biggr).
\end{equation}

%


Step 2. We now use the margin condition to provide lower bounds on
$\underline{\Delta}_j$ for each $j \in\cJ$.
Assume without loss of generality that the indexing of the bins is such
that $\cJ=\{1,\ldots,j_1\}$ and that the gaps are ordered
$0<\underline{\Delta}_1 \le\underline{\Delta}_2 \le\cdots\le
\underline{\Delta}_{j_1}$.
For any $j \in\cJ$, from the definition of $\underline{\Delta}_j$,
there exists a suboptimal arm $i \in\cI_j^0$ such that $\underline
{\Delta}_j=\bar f^\star_j-\bar f_j^{(i)}$. But since the function
$f^\star-f^{(i)}$ satisfies the smoothness condition with parameters
$(\beta,2L)$, we find that if $\underline{\Delta}_j \le\delta$ for
some $\delta>0$, then
\[
0<f^\star(x) -f^{(i)}(x) \le\delta+ 2Ld^{\beta/2}M^{-\beta}\qquad
\forall x \in B_j.
\]
Together with the fact that $f^\star-f^\sharp>0$ over $B_j$ for any
$j \in\cJ$ (see step 1 above), it yields
\[
P_X \bigl[ 0 < f^{\star}-f^{\sharp} \le\underline{
\Delta}_j+2Ld^{\beta/2}M^{-\beta} \bigr]\ge\sum
_{k=1}^{j_1} p_k \1 (0<\underline{
\Delta}_k \le\underline{\Delta}_{j})\ge\frac
{\underline{c} j}{M^d},
\]
where we used the fact that $p_k=P_X(B_k) \ge\underline{c}/M^d$.
Define $j_2 \in\cJ$ to be the largest integer such that $\underline
{\Delta}_{j_2} \le\delta_0/c_1$. Since for any $j \in\cJ$, we have
$\underline{\Delta}_j>M^{-\beta}$, the margin condition yields for
any $ j \in\{1,\ldots, j_2\}$ that
\[
P_X \bigl[ 0 < f^{\star}-f^{\sharp} \le\underline{
\Delta}_j+2Ld^{\beta/2}M^{-\beta} \bigr]\le
C_\delta(c_1\underline{\Delta}_j)^\alpha,
\]
where we have used the fact that $\underline{\Delta}_j+2Ld^{\beta
/2}M^{-\beta}\le
c_1\underline{\Delta}_j \le\delta_0$, for any $j \in\{1,\ldots,
j_2\}$. The previous two inequalities, together with the fact that
$\underline{\Delta}_j > M^{-\beta}$ for any $j \in\cJ$, yield
\[
\underline{\Delta}_{j}\ge c_4 \biggl(
\frac{j}{M^d} \biggr)^{1/\alpha} \vee M^{-\beta}=:
\gamma_j\qquad \forall j \in\{1,\ldots, j_2\}.
\]
Therefore, using the fact that $\underline{\Delta}_j \geq\delta
_0/c_1$ for $j \geq j_2$, we get from (\ref{EQregJ}) that
%
%
\begin{eqnarray}
\label{eqwellbehaved}\quad
&&\sum_{j\in\cJ}\E
\sfR_j(\bpi) \nonumber\\[-8pt]\\[-8pt]
&&\qquad\le c_5 \Biggl[\sum
_{j=1}^{j_2} K \frac{\blog({n}\gamma
_j^2/{M^d} )}{\gamma_j} + \sum
_{j=j_2+1}^{j_1}K\log(n)+ \sum
_{j \in
\cJ} nM^{-d-\beta}q_j \Biggr].\nonumber
\end{eqnarray}

\textit{Fourth part}: \textit{Putting things together.\quad}
Combining (\ref{eqstronglyill}), (\ref{eqweaklyill}) and (\ref
{eqwellbehaved}), we obtain
the following bound:
%
%
\begin{eqnarray}
\label{EQreg2} \E R_n(\bpi)&\le& c_6 \Biggl[
nM^{-\beta(1+\alpha)}\nonumber\\
&&\hspace*{13.5pt}{} + \sqrt{K\log(K)}M^{{d}/{2}-\alpha\beta}\sqrt
{n}+K\sum
_{j=1}^{j_2} \frac
{\blog({n}\gamma_j^2/{M^d} )}{\gamma_j}
\\
&&\qquad\hspace*{72.5pt}{} + K M^d \log n +nM^{-d-\beta}\sum
_{j \in\cJ}q_j \Biggr].\nonumber
\end{eqnarray}

We now bound from above the first sum in (\ref{EQreg2}) by
decomposing it into two terms. From the definition of $\gamma_j$,
there exists an integer $j_3$ satisfying $c_7M^{d-\alpha\beta} \le
j_3 \le2c_7M^{d-\alpha\beta}$ and such that $\gamma_j = M^{-\beta
}$ for $j \leq j_3$ and $\gamma_j = c_4(jM^{-d})^{1/\alpha}$ for $j
>j_3$. It holds
%
%
\begin{equation}
\label{EQintegral1} \sum_{j=1}^{j_3}
\frac{\blog({n}\gamma_j^2/{M^d}
)}{\gamma_j} \leq c_8M^{d+\beta(1-\alpha)}\blog\biggl(
\frac
{n}{M^{2\beta+d}} \biggr)
\end{equation}
and
%
%
\begin{eqnarray}\label{EQintegral}\qquad
\sum_{j=j_3+1}^{j_2} \frac{\blog({n}\gamma
_j^2/{M^d} )}{\gamma_j}
&\leq&
c_9\sum_{j=j_3+1}^{M^d} \biggl(
\frac
{j}{M^d} \biggr)^{-{1}/{\alpha}}\blog\biggl(\frac{n}{M^d} \biggl[
\frac{j}{M^d} \biggr]^{{2}/{\alpha}} \biggr)
\nonumber\\[-8pt]\\[-8pt]
&\le& c_{10} M^d \int_{M^{-\alpha\beta}}^{1}
\blog\biggl(\frac{n}{M^d}x^{{2}/{\alpha}} \biggr)x^{-1/\alpha} \,\ud
x. \nonumber
\end{eqnarray}
Since $\alpha< 1$, this integral is bounded by $c_{10}M^{\beta
(1-\alpha)} (1+\blog(n/M^{2\beta+d} ) ) $.

The second sum in (\ref{EQreg2}) can be bounded as
%
%
\begin{eqnarray}\label{EQsumqj}
\sum_{j \in\cJ} q_j &=& \sum
_{j \in\cJ} \p\bigl\{0< f^{\star
}(X)-f^{\sharp}(X)
\leq c_1M^{-\beta} | X \in B_j \bigr\}
\nonumber\\[-8pt]\\[-8pt]
&\leq&\frac{M^d}{\underline{c}}\p\bigl\{0<
f^{\star}(X)-f^{\sharp
}(X)\leq c_1M^{-\beta} \bigr
\} \leq\frac{c_1^\alpha}{\underline
{c}} M^ {d-\beta\alpha}.\nonumber
\end{eqnarray}
Putting together (\ref{EQreg2})--(\ref{EQsumqj}), we obtain
\begin{eqnarray*}
\E R_n(\bpi)&\le& c_{11} \biggl[nM^{-\beta(1+\alpha)}+
\sqrt{K\log(K)}M^{{d}/{2}-\alpha\beta}\sqrt{n}+ KM^{d+\beta
(1-\alpha)}
\\
&&\hspace*{67.6pt}{} +KM^{d+\beta(1-\alpha)}\blog\biggl(\frac
{n}{M^{2\beta+d}} \biggr)+KM^d
\log n \biggr],
\end{eqnarray*}
and the result follows by choosing $M$ as prescribed.
\end{pf*}

We should point out that the version of the \textsc{bse} described
above specifies the number of bins $M$ as a function of the
horizon $n$, while in practice one may not have foreknowledge of
this value. This limitation can be easily circumvented by using
the so-called \textit{doubling argument} (see, e.g., page 17 in~\cite
{CesLug06}) which
consists of ``reseting'' the game at times $2^{k}, k=1, 2, \ldots\,$.

The reader will note that when $\alpha=1$ there is a potentially
superfluous $\log n$ factor appearing in the upper bound using the same
proof. More generally, for any $\alpha\ge1$, it is possible to
minimize the expression in (\ref{EQreg2}) with respect to $M$, but
the optimal value
of $M$ would then depend on the value of $\alpha$. This sheds
some light on a significant limitation of the \textsc{bse} which
surfaces in this parameter regime: for $n$ large enough, it requires
the operator to
pull each arm at least once in each bin and therefore to incur an expected
regret of at least order $M^d$. In other words, the \textsc{bse}
splits the space $\cX$ in ``too many'' bins when $\alpha\ge1$.
Intuitively this can be understood as follows. When $\alpha\ge1$,
the gap function $f^\star(x)-f^\sharp(x)$ is bounded away from zero
on a large subset of $\cX$. Hence there is no need to carefully
estimate it since the optimal arm is the same across the region. As a
result, one could use larger bins in such regions
reducing the overall number of bins and therefore removing the
extra logarithmic term alluded to above.

%

\section{Adaptively Binned Successive Elimination}\label{SECabse}


We need the following definitions. Assume that $n \ge K\log(K)$ and
let $k_0$ be the smallest integer such that
%
%
\begin{equation}
\label{EQdefk0} 2^{-k_0} \le\biggl(\frac{K\log(K)}{n}
\biggr)^{{1}/({d+2\beta})}.
\end{equation}
For any bin $B \in\bigcup_{k=0}^{k_0} \cB_{2^k}$, let $\ell_B$ be
the smallest integer such that
%
%
\begin{equation}
\label{EQlbeps} U\bigl(\ell_B, n|B|^{d}\bigr)
\le2c_0|B|^{\beta},
\end{equation}
where $U$ is defined in (\ref{EQdefU}) and $c_0=2Ld^{\beta/2}$. This
definition implies that
%
%
\begin{equation}
\label{EQboundlB} \ell_B \le C_\ell|B|^{-2\beta}
\log\bigl(n|B|^{(2\beta+d)}\bigr)
\end{equation}
for $C_\ell>0$, because $x \mapsto U(x,n|B|^d)$ is decreasing for $x>0$.

The \textsc{abse} policy operates akin to the \textsc{bse} except that
instead of fixing a partition~$\cB_M$, it relies on an adaptive
partition that is refined over time. This partition is better
understood using the notion of rooted tree.

Let $\cT^*$ be a tree with root $\cX$ and maximum depth $k_0$. A node
$B$ of $\cT^*$ with depth $k=0,\ldots, k_0-1$ is a set from the
regular partition $\cB_{2^k}$. The children of node $B \in\cB_{2^k}$
are given by $\burst(B)$, defined to be the collection of $2^d$ bins
in $\cB_{2^{k+1}}$ that forms a partition of $B$.

Note that the set $\cL$ of leaves of each subtree $\cT$ of $\cT^*$
forms a partition of $\cX$.
The \textsc{abse} policy constructs a sequence of partitions $\cL
_1,\ldots, \cL_n$ that are leaves of subtrees of $\cT^*$. At a given
time $t=1,\ldots, n$, we refer to the elements of the current
partition $\cL_t$ as \textit{live} bins. The sequence of partitions is
nested in the sense that if $B \in\cL_t$, then either $B \in\cL
_{t+1}$ or $\burst(B) \subset\cL_{t+1}$. The sequence $\cL_1,\ldots, \cL
_n$ is constructed as follows.

In the initialization step, set $\cL_0=\varnothing$, $\cL_1=\cX$, and
the initial set of arms $\cI_{\cX}=\{1, \ldots, K\}$. Let $t\le n$ be
a time such that $\cL_t \neq\cL_{t-1}$, and let $\mathsf{B}_t$ be
the collection of sets $B$ such that $B \in\cL_t \setminus\cL
_{t-1}$. We say that the bins $B \in\mathsf{B}_t$ are \textit{born} at
time $t$. For each set $ B \in\mathsf{B}_t$, assume that we are given
a set of active arms $\cI_{B}$. Note that $t=1$ is such a time with
$\mathsf{B}_1=\{\cX\}$ and active arms $\cI_\cX$. For each born bin
$B \in\mathsf{B}_t$, we run a \textsc{se} policy $\hpi_B$ initialized
at time $t$ with initial set of arms $\cI_B$ and parameters
$T_B=n|B|^{-d}$, $\gamma=2$. Such a policy is defined in Section \ref
{subbin}. Let $t(B)$ denote the time at which $\hpi_B$ has reached
$\ell_B$ rounds and let
%
%
\begin{equation}
\label{EQtildeSB} \tilde N_B(t)=\sum
_{l=1}^t \1(X_t \in B, B \in
\cL_t)
\end{equation}
denote the number of times covariate $X_t$ fell in bin $B$ while $B$
was a live $B$. At time $t(B)+1$, we replace the node $B$ by its
children $\burst(B)$ in the current partition. Namely, $\cL
_{t(B)+1}=(\cL_{t(B)}\setminus B)\cup\burst(B)$. Moreover, to each
bin $B' \in\burst(B)$, we assign the set $\cI_{B'}=\S_{B, \tilde
N_{B}(t(B))}$ of arms that were left active by policy $\hpi_B$ on its
parent $B$ at the end of the $\ell_B$ rounds. This procedure is
repeated until the horizon $n$ is reached.

The intuition behind this policy is the following. The parameters of
the \textsc{se} policy $\hpi_B$ run at the birth of bin $B$ are chosen
exactly such that arms $i$ with average gap $|\bar f^\star_B -\bar
f^{(i)}_B| \ge C|B|^\beta$ are eliminated by the end of $\ell_B$ rounds
with high probability. The smoothness condition ensures that these
eliminated arms satisfy $f^\star(x) >f^{(i)}(x)$ for all $x \in B$ so
that such arms are uniformly suboptimal on bin $B$. Among the kept
arms, none is uniformly better than another, so bin $B$ is burst and
the process is repeated on the children of $B$ where other arms may be
uniformly suboptimal. The formal definition of the \textsc{abse} is given
in Policy~\ref{ALGabse}; it satisfies the following theorem.

\begin{algorithm}[t]
\caption{Adaptively Binned Successive Elimination (\textsc{abse})}
\label{ALGabse}
\begin{algorithmic}
\REQUIRE Set of arms $\cI_{\cX}=\{1,\ldots, K\}$. Parameters $n,
c_0=2Ld^{\beta/2}, k_0$.
\ENSURE$\tpi_1,\ldots, \tpi_n \in\cI$.
\STATE$t \leftarrow0$, $k \leftarrow0$, $\cL\leftarrow\{\cX\}$.
\STATE Initialize a \textsc{se} policy $\hpi_\cX$ with parameters $T=n,
\gamma=2$ and arms $\cI=\cI_\cX$.
\STATE$N_\cX\leftarrow0$.
\FOR{$t=1,\ldots, n$}
\STATE$B \leftarrow\cL(X_t)$.
\STATE$N_{B} \leftarrow N_B + 1$. \hfill\textttt{/count times $X_t
\in B$/}
\STATE$\tpi_t \leftarrow\hpi_{B,N_B}$ (observe $Y^{(\tpi_t)}_t$).
\hfill\textttt{ /choose arm from \textsc{se} policy $\hpi_B$/}
\STATE$\tau_B \leftarrow\hat\tau_{B,N_B}$ \hfill\textttt{/update
number of rounds for $\hpi_B$/}
\STATE$\cI_B \leftarrow\S_{B,N_B}$ \hfill\textttt{/update active
arms for $\hpi_B$/}
\IF{$\tau_B \ge\ell_B$ and $|B|\ge2^{-k_0+1}$ and $|\cI_B|\ge2$
\hfill\textttt{/conditions to $\burst(B)$/}}
\FOR{$B' \in\burst(B)$}
\STATE$\cI_{B'} \leftarrow\cI_{B}$\hfill\textttt{/assign remaining
arms as initial arms/}
\STATE Initialize \textsc{se} policy $\hpi_{B'}$ with $T=n|B'|^d,
\gamma=2$ and arms $\cI=\cI_{B'}$.
\STATE$N_{B'}\leftarrow0$.\hfill\textttt{/set time to 0 for new \textsc
{se} policy/}
\ENDFOR
\STATE$\cL\leftarrow\cL\setminus B$ \hfill\textttt{/remove $B$
from current partition/}
\STATE$\cL\leftarrow\cL\cup\burst(B)$\hfill\textttt{/add $B$'s
children to current partition/}
\ENDIF
\ENDFOR
\end{algorithmic}
\end{algorithm}

%
\begin{theorem}
\label{THmainadap} Fix $\beta\in(0,1]$, $L>0$, $\alpha>0$, assume that
$n \ge K\log(K)$ and consider a problem in
$\cM^K_{\cX}(\alpha,\beta,L)$. If $\alpha<\infty$, then the \textsc{abse}
policy $\tpi$ has an expected regret at time $n$ bounded by
\[
\E R_n(\tpi)\le C n \biggl(\frac{K\log(K)}{n} \biggr)^{{\beta
(\alpha+1)}/({2\beta+d})},
\]
where $C>0$ does not depend on $K$. 
If $\alpha=\infty$, then $\E R_n(\tpi)\le C K\log(n)$.
\end{theorem}
Note that the bounds given in Theorem~\ref{THmainadap} are optimal in
a minimax sense when $K=2$. Indeed, the lower bounds of \cite
{AudTsy07} and~\cite{RigZee10} imply that the bound on expected regret
cannot be improved as a function of $n$ except for a constant
multiplicative term. The lower bound proved in~\cite{AudTsy07} implies
that any policy that received information from \textit{both} arms at
each round has a regret bound at least as large as the one from
Theorem~\ref{THmainadap}, up to a multiplicative constant. As a
result, there is no price to pay for being in a partial information
setup and one could say that the problem of nonparametric estimation
dominates the problem associated to making decisions sequentially.

Note also that when $\alpha=\infty$, Proposition \ref
{propalpha-beta} implies that there exists a unique optimal arm over
$\cX$ and that all other arms have reward bounded away from that of
the optimal arm. As a result, given this information, one could operate
as if the problem was static by simply discarding the covariates.
Theorem~\ref{THmainadap} implies that in this case, one recovers the
traditional regret bound of the static case without the knowledge that
$\alpha=\infty$.

\begin{pf*}{Proof of Theorem~\ref{THmainadap}}
We first consider the case where $\alpha< \infty$, which
implies that $\alpha\beta\leq d$; see Proposition~\ref{propalpha-beta}.

We keep track of positive constants by numbering them $c_1,c_2,\ldots\,$,
yet they might differ from previous sections.
On each newly created bin $B$, a new \textsc{se} policy is initialized,
and we denote by $Y_{B, 1}^{(i)}, Y_{B, 2}^{(i)}, \ldots\,$, the rewards
obtained by successive pulls of a remaining arm $i$. Their average
after $\tau$ rounds/pulls is denoted by
\[
\bar Y_{B, \tau}^{(i)}:=\frac{1}{\tau}\sum
_{s=1}^\tau Y_{B,
s}^{(i)}.
\]
For any 
integer $s$, define $\eps_{B,s}=2U(s, n|B|^d)$, where $U$ is defined
in (\ref{EQdefU}).

For any $B \in\mathcal{T}^*\setminus\{\cX\}$, define the unique
\textit
{parent} of $B$ by
\[
\pp(B):=\bigl\{B'\in\mathcal{T}^*\dvtx B \in\burst
\bigl(B'\bigr)\bigr\}
\]
and $\pp(\cX)=\varnothing$.
Moreover, let $\pp^1(B)=\pp(B)$ and for any $k \ge2$ define
recursively $\pp^k(B)=\pp(\pp^{k-1}(B))$. Then the set of \textit
{ancestors} of any $B \in\mathcal{T}^*$ is denoted by $\cP(B)$ and
defined by
\[
\cP(B)=\bigl\{B'\in\mathcal{T}^*\dvtx B'=
\pp^k(B) \mbox{ for some } k \ge1\bigr\}.
\]
%
Denote by $r^{\mathrm{live}}_n(B)$ the regret incurred by the \textsc
{abse} policy $\tpi$ when covariate $X_t$ fell in a \textit{live} bin $B
\in\cL_t$, where we recall that $\cL_t$ denotes the current
partition at time $t$.
It is defined by
\[
r^{\mathrm{live}}_n(B)=\sum_{t=1}^n
\bigl[f^\star(X_t) - f^{(\tpi
_t(X_t))}(X_t)\bigr]
\1 (X_t \in B)\1(B\in\cL_t).
\]
We also define $\cB_t:= \bigcup_{s \leq t} \cL_s$ to be the set of
bins that were born at some time $s \leq t$. We denote by $r_n^{\mathrm
{born}}(B)$ the regret incurred when covariate $X_t$ fell in such a
bin. It is defined by
\[
r^{\mathrm{born}}_n(B)=\sum_{t=1}^n
\bigl[f^\star(X_t) - f^{(\tpi
_t(X_t))}(X_t)\bigr]
\1 (X_t \in B)\1(B \in\cB_t).
\]
Observe that if we define $\tilde r_n:= r_n^{\mathrm{born}}(\cX)$, we
have $\E R_n(\tpi)=\E\tilde r_n$ since $\cX\in\cB_t$ and $X_t \in
\cX$ for all $t$.
Note that for any $B \in\mathcal{T}^*$,
%
%
\begin{equation}
\label{EQdecompborn} r^{\mathrm{born}}_n(B)=r^{\mathrm{live}}_n(B)+
\sum_{B' \in\burst
(B)} r^{\mathrm{born}}_n
\bigl(B'\bigr).
\end{equation}
Denote by $\cI_B=\S_{B,t_B}$ the set of arms left active by the \textsc
{se} policy $\hpi_B$ on $B$ at the end of $\ell_B$ rounds. Moreover,
define the following reference sets of arms:
\[
\ucI:= \Bigl\{i \in\{1,\ldots, K\}\dvtx \sup_{x \in B }
f^\star(x)- f^{(i)}(x) \le c_0 |B|^{\beta}
\Bigr\},
\]
%
%
\[
\ocI:= \Bigl\{i \in\{1,\ldots, K\}\dvtx \sup_{x \in B }
f^\star(x)- f^{(i)}(x) \le8c_0 |B|^{\beta}
\Bigr\}.
\]
Define the event $\cA_B:=\{ \ucI\subseteq\cI_B \subseteq\ocI\}$
on which the remaining arms have a gap of the correct order and observe
that (\ref{EQdecompborn}) implies that
\[
r^{\mathrm{born}}_n(B)=r^{\mathrm{born}}_n(B)\1\bigl(
\cA_B^c\bigr)+r^{\mathrm
{live}}_n(B)\1(
\cA_B)+\sum_{B' \in\burst(B)} r^{\mathrm
{born}}_n
\bigl(B'\bigr)\1( \cA_{B}).
\]
Let $\mathcal{L}^*$ denote the set of leaves of $\mathcal{T}^*$, that
is the set of
bins $B$ such that $|B|=2^{-k_0}$. In what follows, we adapt the
convention that $\prod_{B' \in\cP(\cX)} \1(\cA_{B'})=1$.

We are going to treat regret incurred on live nonterminal nodes and
live leaves separately and differently. As a result, the quantity we
are interested in is decomposed as $\tilde r_n=\tilde r_n(\mathcal{T}^*
\setminus\mathcal{L}^*) + \tilde r_n(\mathcal{L}^*)$ where
\[
\tilde r_n\bigl(\mathcal{T}^*\setminus\mathcal{L}^*\bigr):=\sum
_{B \in
\mathcal{T}^*\setminus\mathcal{L}^*} \bigl(r^{\mathrm{born}}_n(B) \1
\bigl(\cA_B^c\bigr)+r^{\mathrm{live}}_n(B) \1
(\cA_B) \bigr)\prod_{B' \in\cP(B)} \1(
\cA_{B'})
\]
is the regret accumulated on live nonterminal nodes, and
\[
\tilde r_n\bigl(\mathcal{L}^*\bigr):=\sum
_{B \in\mathcal{L}^*} r^{\mathrm
{born}}_n(B)\prod
_{{B' \in\cP(B)}}\1(\cA_{B'})=\sum
_{B \in\mathcal{L}^*} r^{\mathrm
{live}}_n(B)\prod
_{{B' \in\cP(B)}}\1(\cA_{B'})
\]
is regret accumulated on live leaves. Our proof relies on the following
events: $\cG_B:=\bigcap_{B' \in\cP(B)}\cA_{B'}$.\vspace*{9pt}

\textit{First part}: \textit{Control of the regret on the
nonterminal nodes.\quad}
Fix \mbox{$B \in\mathcal{T}^*\setminus\mathcal{L}^*$}. On $\cG_B$, we
have $\cI_{\pp(B)}
\subseteq\bar{\cI}_{\pp(B)}$ so that any active arm $i \in\cI
_{\pp(B)}$ satisfies $\sup_{x \in\pp(B) } |f^\star(x)- f^{(i)}(x)|
\le8c_0 |\pp(B)|^{\beta}$. Moreover, regret is only incurred at
points where $f^*-f^\sharp>0$, so\vadjust{\goodbreak} defining $c_1:=2^{3+\beta}c_0$ and
conditioning on events $\{X_t \in B\}$ yields
\[
\E\bigl[r^{\mathrm{live}}_n(B) \1(\cG_B\cap
\cA_B) \bigr]\le\E\bigl[\tilde N_B(n)
\bigr]c_1|B|^\beta q_B \le c_1 K
\ell_B |B|^\beta q_B,
\]
where $q_B=P_X (0<f^\star-f^\sharp\le c_1 |B|^\beta| X
\in B )$ and $\tilde N_B(n)$ is defined in (\ref{EQtildeSB}).

We can always assume that $n$ is greater than $n_0 \in\N$, defined by
\[
n_0= \biggl\lceil K\log(K) \biggl(\frac{c_1}{\delta_0}
\biggr)^{
({d+2\beta})/{\beta}} \biggr\rceil\qquad\mbox{so that } c_12^{-k_0\beta}
\leq\delta_0,
\]
and let $k_1 \leq k_0$ be the smallest integer such that
$c_12^{-k_1\beta} \leq\delta_0$. Indeed, if $n \leq n_0$, the result
is true with a constant large enough.

Applying the same argument as in (\ref{EQsumqj}) yields the existence
of $c_2>0$ such that, for any $k \in\{0,\ldots,k_0\}$,
\[
\sum_{|B|=2^{-k}} q_B \leq
c_22^{k(d-\beta\alpha)}.
\]
Indeed, for $k \geq k_1$ one can define $c_2= c_1^\alpha/\underline
{c}$, and the same equation holds with $c_2=2^{dk_1}$ if $k \leq k_1$.
Summing over all depths $k \leq k_0-1$, we obtain
%
%
\begin{eqnarray}
\label{EQhatrn1}
&&
\E\biggl[ \sum_{B \in\mathcal{T}^*\setminus\mathcal{L}^*}
r^{\mathrm
{live}}_n(B) \1(\cG_B\cap
\cA_B) \biggr]\nonumber\\[-10pt]\\[-10pt]
&&\qquad\leq c_1c_2C_\ell K
\sum_{k=0}^{k_0-1} 2^{k(d+\beta-\alpha\beta)}\log
\bigl(n2^{-k(2\beta
+d)} \bigr).\nonumber
\end{eqnarray}

On the other hand, for every bin $B \in\mathcal{T}^*\setminus
\mathcal{L}^*$, one also has
%
%
\begin{equation}
\label{EQrn1} \E\bigl[r^{\mathrm{born}}_n(B) \1\bigl(
\cG_B\cap\cA_B^c\bigr) \bigr]\le
c_1 n |B|^\beta q_BP_X(B) \p
\bigl(\cG_B\cap\cA_B^c\bigr).
\end{equation}

It remains to control the probability of $\cG_B \cap\cA_B^c$; we
define $\p^{\cG_B}(\cdot):=\p( \cdot\cap\cG_B)$.
On $\cG_B$, the event $ \cA_{B}^c$ can occur in two ways:
\begin{longlist}[(ii)]
\item[(i)]By eliminating an arm $i \in\ucI$ at the end of the at
most $\ell_B$ rounds played on bin $B$. These arms satisfy $\sup_{x
\in B}f^\star(x)-f^{(i)}(x) < c_0|B|^{\beta}$; this event is denoted
by $\cD_B^1$.
\item[(ii)]By not eliminating an arm $i \notin\ocI$ within the at
most $\ell_B$ rounds played on bin $B$. These arms satisfy $\sup_{x
\in B}f^\star(x)-f^{(i)}(x) \ge8c_0|B|^{\beta}$; this event is
denoted by $\cD_B^2$.
\end{longlist}
We use the following decomposition:
%
%
\begin{equation}
\label{EQdecompG} \p^{\cG_B}\bigl(\cA_{B}^c
\bigr)= \p^{\cG_B}\bigl(\cD_B^1\bigr) +
\p^{\cG_B}\bigl(\cD_B^2 \cap\bigl(
\cD_B^1\bigr)^c\bigr).
\end{equation}

We first control the probability of making error (i). Note that for any
$s \le\ell_B$ and any arms $i \in\ucI, i' \in\cI_{\pp(B)}$, it holds
\[
\bar f^{(i')}_B - \bar f^{(i)}_B
\le\bar f^\star_B- \bar f^{(i)}_B <
c_0|B|^{\beta} \le\frac{\varepsilon_{B,\ell_B}}{2}.
\]
Therefore, if an arm $i \in\ucI$ is eliminated, that is, if there
exists $ i' \in\cI_{\pp(B)}$ such that $ \bar{Y}_{B,s}^{(i')}-\bar
{Y}_{B,s}^{(i)}>\varepsilon_{B,s}$ for some $s \le\ell_B$, then
either $\bar{f}^{(i)}_B$ or $\bar{f}^{(i')}_B$ does not\vadjust{\goodbreak} belong to its
respective confidence interval $ [\bar{Y}_{B, s}^{(i)} \pm
\varepsilon_{B,s}/4 ]$ or $ [\bar{Y}_{B, s}^{(i')} \pm
\varepsilon_{B,s}/4 ]$ for some $s \le\ell_B$.
%
Therefore, since $-\bar{f}^{(i)}_B\leq Y_s -\bar{f}^{(i)}_B \leq
1-\bar{f}^{(i)}_B$,
%
%
\begin{equation}
\label{EQboundG1}\qquad \p^{\cG_B}\bigl(\cD_B^1
\bigr) \le\p\biggl\{\exists s \leq\ell_B; \exists i \in
\cI_{\pp(B)}; \bigl\llvert\bar{Y}_s^{(i)}-
\bar{f}^{(i)}_B\bigr\rrvert\geq\frac{\varepsilon_{B,s}}{4} \biggr\}
\leq2K\frac{\ell
_B}{n|B|^d},
\end{equation}
where in the second inequality, we used Lemma~\ref{LEMpeeling}.
%

Next, we\vspace*{1pt} treat error (ii). For any $i \notin\ocI$, there exists $x^{(i)}
$ such that $f^\star(x^{(i)})-f^{(i)}(x^{(i)}) > 8c_0|B|^\beta$. Let
$\check{\imath} =\check{\imath}(i) \in\cI$ be any arm such that
$f^\star(x^{(i)})=f^{(\check{\imath})}(x^{(i)})$; the smoothness
condition implies that
%
%
\begin{eqnarray}
\label{EQgapG2} \bar f^{(\check{\imath})}_B &\ge& f^{(\check{\imath})}
\bigl(x^{(i)}\bigr) -c_0|B|^{\beta} > f^{(i)}
\bigl(x^{(i)}\bigr) + 7c_0|B|^\beta\nonumber\\[-9pt]\\[-9pt]
&\ge&\bar
f^{(i)}_B + 6c_0|B|^\beta\geq\bar
f^{(i)}_B + \tfrac{3}{2} \varepsilon_{B,\ell
_B}.\nonumber
\end{eqnarray}
On the event $(\cD_B^1)^c$, no arm in $\ucI$, and in particular any
of the arms $\check{\imath}(i), i \in\cI_{\pp(B)}\setminus\ocI$,
has been eliminated until round $\ell_B$. Therefore, the event $\cD
_B^2\cap(\cD_B^1)^c$ occurs if there exists $ i \notin\ocI$ such
that $\bar{Y}_{B, \ell_B}^{(\check{\imath})}-\bar{Y}_{B, \ell
_B}^{(i)}\le\varepsilon_{B,\ell_B}$.
In view of (\ref{EQgapG2}) and (\ref{EQlbeps}), it implies that
there exists $i \in\cI_{\pp(B)}$ such that
\[
\bigl|\bar{Y}_{B, \ell_B}^{(i)} - \bar f^{(i)}_B
\bigr|\ge\frac{\varepsilon
_{B,\ell_B}}{4}.
\]
Hence, the probability of error (ii) can be bounded by
%
%
\begin{eqnarray}
\label{EQboundG2} \p^{\cG_B}\bigl(\cD_B^2
\cap\bigl(\cD_B^1\bigr)^c\bigr) &\le&\p\biggl
\{\exists i \in\cI_{\pp(B)}\dvtx \bigl|\bar{Y}_{B, \ell_B}^{(i)} -
\bar f^{(i)}_B \bigr|\ge\frac{\varepsilon_{B,\ell_B}}{4} \biggr\}\nonumber\\[-10pt]\\[-10pt]
&\leq&2K
\frac{\ell
_B}{n|B|^d},\nonumber
\end{eqnarray}
where the second inequality follows from (\ref{EQhoeffding}).

Putting together (\ref{EQdecompG}), (\ref{EQboundG1}), (\ref
{EQboundG2}) and (\ref{EQboundlB}), we get
\[
\p^{\cG_B}\bigl(\cA_{B}^c\bigr) \le4K
\frac{\ell_B}{n|B|^d} \le4C_\ell\frac{K}{n}|B|^{-(2\beta+d)}
\log\bigl(n|B|^{(2\beta+d)}\bigr).
\]
Together with (\ref{EQrn1}), it yields for $B \in\mathcal
{T}^*\setminus\mathcal{L}^*
$ that
\[
\E\bigl[r^{\mathrm{born}}_n(B) \1\bigl(\cG_B\cap
\cA_B^c\bigr) \bigr]\le c_3 K
|B|^{-(\beta+d)}\log\bigl(n|B|^{(2\beta+d)}\bigr) q_BP_X(B).
\]
If $k$ is such that $c_1 2^{-k\beta} > \delta_0$, then any bin $B$
such that $|B|=2^{-k}$ satisfies $
\E[r^{\mathrm{born}}_n(B) \1(\cG_B\cap\cA_B^c) ]\le
c_4 K\log n$.
If $k$ is such that $c_1 2^{-k\beta} \le\delta_0$, then the above
display together with the margin condition yield
\[
\E\biggl[\sum_{|B|=2^{-k}}r^{\mathrm{born}}_n(B)
\1\bigl(\cG_B\cap\cA_B^c\bigr) \biggr]
\le c_5 K 2^{k(\beta+d-\alpha\beta)}\log\bigl(n2^{-k(2\beta+d)}\bigr).
\]
Summing over all depths $k=0,\ldots, k_0-1$ and using (\ref
{EQhatrn1}), we obtain
%
%
\begin{equation}
\label{EQregnodes} \E\bigl[\tilde r_n\bigl(\mathcal{T}^*\setminus
\mathcal{L}^*\bigr)\bigr]\le c_6 K \sum
_{k=0}^{k_0-1}2^{k(\beta+d-\alpha\beta)}\log\bigl(n2^{-k(2\beta+d)}
\bigr).\vadjust{\goodbreak}
\end{equation}
We now compute an upper bound on the right-hand side of the above
inequality. Fix $k=0, \ldots, k_0$ and define
\[
S_{k}=\sum_{j=0}^{k}
2^{j(d+\beta-\beta\alpha)}=\frac
{2^{(k+1)(d+\beta-\beta\alpha)}-1}{2^{d+\beta-\beta\alpha}-1}.
\]
Observe that
\[
2^{k(d+\beta-\beta\alpha)}\log\bigl(n 2^{-k(d+2\beta)} \bigr)=(S_{k}-S_{k-1})
\log\bigl(n[c_{7} S_{k}+1]^{-({d+2\beta
})/({d+\beta-\beta\alpha})} \bigr),
\]
where $c_7:=2^{d+\beta-\beta\alpha}-1$.
Therefore, (\ref{EQregnodes}) can be rewritten as
%
%
\begin{eqnarray}\label{EQregnodes2}
&&
\E\bigl[\tilde r_n\bigl(\mathcal{T}^*\setminus\mathcal{L}^*\bigr)
\bigr]\nonumber\\
&&\qquad\le c_6K \Biggl[\sum_{
k=1}^{k_0-1}
(S_{k}-S_{k-1})\log\bigl(n[c_{7}
S_{k}+1]^{-
({d+2\beta})/({d+\beta-\beta\alpha})} \bigr) +\log n \Biggr]
\nonumber
\\
&&\qquad\le c_6K \biggl[\int_{0}^{S_{k_0-1}}
\log\bigl(n[c_{7} x+1]^{-
({d+2\beta})/({d+\beta-\beta\alpha})} \bigr) \,\ud x +\log n \biggr]
\\
&&\qquad\le c_8K \bigl[2^{k_0(d+\beta-\beta\alpha)} \log\bigl
(n2^{-k_0(d+2\beta)}
\bigr) + \log n \bigr]
\nonumber
\\
&&\qquad\le c_9 n \biggl(\frac{n}{K\log(K)} \biggr)^{-{\beta
(1+\alpha)}/({d+2\beta})},\nonumber
\end{eqnarray}
where we used (\ref{EQdefk0}) in the last inequality and the fact
that $\log(n)$ is dominated by $n^{1-\beta(1+\alpha)/(d+2\beta)}$
since $\alpha\beta\leq d$.\vspace*{9pt}


%

\textit{Second part}: \textit{Control of the regret on the leaves.\quad}
Recall that the set of leaves $\mathcal{L}^*$ is composed of bins $B$
such that
$|B|=2^{-k_0}$. Proceeding in the same way as in~(\ref{EQrn1}), we
find that for any $B \in\mathcal{L}^*$, it holds
\[
\E\bigl[r^{\mathrm{live}}_n(B) \1(\cG_B) \bigr]\le
c_1 n |B|^\beta P_X\bigl(0<f^\star-f^\sharp
\le c_1 |B|^\beta, X \in B\bigr).
\]
Since $n \geq n_0$, then $c_12^{-k_0\beta}\leq\delta_0$ and using
the margin assumption, we find
%
%
\begin{eqnarray}
\label{EQregleaves} \sum_{B \in\mathcal{L}^*}\E
\bigl[r^{\mathrm{live}}_n(B) \1(\cG_B) \bigr] &\le&
c_1 n 2^{-k_0\beta(1+\alpha)} \nonumber\\[-8pt]\\[-8pt]
&\le& c_1 n \biggl(
\frac{n}{K\log
(K)} \biggr)^{-{\beta(1+\alpha)}/({d+2\beta})},\nonumber
\end{eqnarray}
where we used (\ref{EQdefk0}) in the second inequality.\vspace*{9pt}

The theorem follows by summing (\ref{EQregnodes2}) and (\ref
{EQregleaves}). If $\alpha=+\infty$, then the same proof holds except
that $\log(n)$ dominates $2^{k_0(\beta+d-\alpha\beta)}\log
(n2^{-k_0(2\beta+d)})$ in~(\ref{EQregnodes2}).
\end{pf*}

%
\begin{appendix}\label{app}

\section*{Appendix: Technical lemma}
\label{SECprstatic}

The following lemma is central to our proof of Theorem \ref
{THstatic}. We recall that a process $Z_t$ is a martingale difference
sequence if $\E[Z_{t+1} |Z_1,\ldots,Z_t ]=0$. Moreover,
if $ a\leq Z_t \leq b$ and if we denote the sequence of averages by
$\bar{Z}_t=\frac{1}{t} \sum_{s=1}^t Z_s$, then Hoeffding--Azuma's
inequality yields that, for every integer $T \ge1$,
%
%
\setcounter{equation}{0}
\begin{equation}
\label{EQhoeffding} \p\biggl\{ \bar{Z}_T \geq\sqrt{
\frac{(b-a)^2}{2T} \log\biggl(\frac{1}{\delta} \biggr)} \biggr\} \leq
\delta.
\end{equation}
The following lemma is a generalization of this result:
%
%
\begin{lem}\label{LEMpeeling}
Let $Z_t$ be a martingale difference sequence with $ a\leq Z_t \leq b$
then, for every $\delta>0$ and every integer $T \ge1$,
\[
\p\biggl\{\exists t \leq T, \bar{Z}_t \geq\sqrt{
\frac
{2(b-a)^2}{t} \log\biggl(\frac{4}{\delta}\frac{T}{t} \biggr)}
\biggr\} \leq\delta.
\]
\end{lem}
\begin{pf}
Define $\varepsilon_t=\sqrt{\frac{2(b-a)^2}{t} \log
(\frac{4}{\delta}\frac{T}{t} )}$. Recall first the
Hoeffding--Azuma maximal concentration inequality. For every $\eta>0$
and every integer $t \ge1$,
%
\[
\p\{\exists s \leq t, s\bar{Z}_s \geq\eta\} \leq\exp\biggl(-
\frac{2\eta^2}{t(b-a)^2} \biggr).
\]
Using a peeling argument, one obtains
\begin{eqnarray*}
\p \{\exists t \leq T, \bar{Z}_t \geq\varepsilon_t \}
&\leq&\sum_{m=1}^{\lfloor\log_2(T)\rfloor} \p\Biggl\{\bigcup
_{ t
=2^m}^{ 2^{m+1}-1}\{ \bar{Z}_t \geq
\varepsilon_t\} \Biggr\}
\\
&\leq&
\sum_{m=1}^{\lfloor\log_2(T)\rfloor} \p\Biggl\{\bigcup
_{ t
=2^m}^{ 2^{m+1}}\{ \bar{Z}_t \geq
\varepsilon_{2^{m+1}}\} \Biggr\} \\
&\leq&\sum_{m=1}^{\lfloor\log
_2(T)\rfloor}
\p\Biggl\{\bigcup_{ t =2^m}^{
2^{m+1}}\bigl\{ t
\bar{Z}_t \geq2^m\varepsilon_{2^{m+1}}\bigr\}
\Biggr\}
\\
&\leq&
\sum_{m=1}^{\lfloor\log_2(T)\rfloor} \exp\biggl(-
\frac
{2 (2^m\varepsilon_{2^{m+1}} )^2}{2^{m+1}(b-a)^2} \biggr)\\
&=&\sum
_{m=1}^{\lfloor\log_2(T)\rfloor}
\frac{2^{m+1}}{T}\frac
{\delta}{4}\le\frac{2^{\log_2(T)+2}}{T}\frac{\delta}{4}\le
\delta.
\end{eqnarray*}
Hence the result.
\end{pf}
\end{appendix}



%

\printaddresses

\end{document}